\documentclass[12pt]{article}
\usepackage{amsmath,amssymb}
\usepackage{graphicx}
\newcommand{\eqdef}{\stackrel{\text{def}}{=}}
\newcommand{\eqdefrm}{\stackrel{\text{\rm def}}{=}}
\newcommand{\n}{\nonumber\\}
\newcommand{\bm}{\boldsymbol}
\newcommand{\ignore}[1]{}
\numberwithin{equation}{section}
\newcommand{\Romannumeral}[1]{\uppercase\expandafter{\romannumeral#1}}

\newcommand{\II}{\text{\Romannumeral{2}}}

\newtheorem{theo}{\bf Theorem}[section]

\newcommand{\cP}{\mathcal{P}}
\newcommand{\cX}{\mathcal{X}}

\allowdisplaybreaks[4]

\begin{document}

\baselineskip=20pt

\newcommand{\preprint}{
\vspace*{-20mm}
   \begin{flushright}\normalsize \sf
    DPSU-21-1\\
  \end{flushright}}
\newcommand{\Title}[1]{{\baselineskip=26pt
  \begin{center} \Large \bf #1 \\ \ \\ \end{center}}}
\newcommand{\Author}{\begin{center}
  \large \bf Satoru Odake$^{\,a}$ and Ryu Sasaki$^{\,b}$ \end{center}}
\newcommand{\Address}{\begin{center}
     $^a$ Faculty of Science, Shinshu University,
     Matsumoto 390-8621, Japan\\
     $^b$ Department of Physics, Tokyo University of Science,
     Noda 278-8510, Japan
   \end{center}}
\newcommand{\Accepted}[1]{\begin{center}
  {\large \sf #1}\\ \vspace{1mm}{\small \sf Accepted for Publication}
  \end{center}}

\preprint
\thispagestyle{empty}

\Title{ Markov Chains Generated by Convolutions\\
of Orthogonality Measures}

\Author

\Address
\vspace{1cm}

\begin{abstract}
About two dozens of exactly solvable Markov chains on one-dimensional finite
and semi-infinite integer lattices are constructed in terms of convolutions
of orthogonality measures of the Krawtchouk, Hahn, Meixner, Charlier, $q$-Hahn,
$q$-Meixner and little $q$-Jacobi polynomials.
By construction, the stationary probability
distributions, the complete sets of eigenvalues and eigenvectors are provided
by the polynomials and the orthogonality measures. An interesting property
possessed by these stationary probability distributions, called `convolutional
self-similarity,' is demonstrated.
\end{abstract}

\section{Introduction}
\label{sec:intro}

In this paper we study one mathematical aspect of Markov chains on
one-dimensional integer lattices. The goal is constructing plenty of examples
of workable Markov chains containing enough adjustable parameters so that rich
and functional applications could be made in diverse disciplines.
Markov chains are most easy to handle when, on top of the stationary
probability distributions, the complete set of eigenvectors and corresponding
eigenvalues are known.
The eigenvectors discussed in the present paper form an orthogonal basis 
as they belong to real symmetric matrices \eqref{realsym}.

We pursue this goal within the framework of orthogonal polynomials of a
discrete variable \cite{nikiforov}--\cite{gasper}. To be more specific, we
construct the basic transition matrix $K$ by certain convolutions of the
`orthogonality measures $=$ stationary probability distributions.'
The information of the stationary distributions and the orthogonal eigenvectors
are built in the scheme. In most cases, the eigenvalues can be extracted
directly from the convolutions. In all cases, the eigenvalues are calculated
exactly, see {\bf Theorem\,\ref{theo2}}. In some respects, in particular, in
incorporating many adjustable parameters, the present method is more
advantageous than some sophisticated procedures, for example those based on the
finite group actions. The presentation of this paper is simple and plain so
that non-experts can understand.

This paper is organised as follows.
In section two, after a common problem setting, the basic properties of
$K(x,y)$, the transition probability matrix from $y$ to $x$, are stated
in {\bf Lemma} and four {\bf Theorems}.
{\bf Lemma} states $K(x,y)$ is triangular in a certain basis.
The most fundamental one is {\bf Theorem\,\ref{theo1}} stating $K$ is
related to a real symmetric matrix $\mathcal{H}$ by a similarity transformation.
The complete set of eigenvectors of $K$ is identified based on Lemma.
It is followed by the eigenvalue formula {\bf Theorem\,\ref{theo2}}, the
spectral representation {\bf Theorem\,\ref{theo3}} and solutions of the initial
value problem and the $\ell$-step transition matrix in
{\bf Theorem\,\ref{theo4}}.
Section three provides fundamental data. The list of five types of convolutions
can be found
in \S\,\ref{sec:convlist}, the polynomial data in \S\,\ref{sec:data}, the
convolutional self-similarity of various stationary distributions are explored
in \S\,\ref{sec:selfsim}.
Many explicit examples of $K(x,y)$ constructed by convolutions are demonstrated
in \S\,\ref{convtype1}--\S\,\ref{convtype5}, corresponding to each type of
convolutions. Markov chains constructed by convolutions of type
(\romannumeral1) and (\romannumeral3) have the special property that the
eigenvalues are directly obtained from the determinant of $K$'s, which have a
factorised form of an upper and lower triangular matrix. The polynomials
defined on a one-dimensional integer lattice $\cX=\{0,1,\ldots,N\}$, the
Krawtchouk, Hahn and $q$-Hahn are the main players. Those $K$'s lead to a
wider class of Markov chains on a semi-infinite lattice by limiting procedures
$N\to\infty$, described by ($q$-)Meixner and Charlier polynomials.
Markov chains constructed by repeated convolutions of type (\romannumeral1) and
(\romannumeral3) for the Krawtchouk are presented in \S\,\ref{sec:(i)(iii)gen}.
Several examples of one parameter families of commuting $K$'s are presented in
\S\,\ref{sec:powerL}.
Section five deals with two related topics. The dual convolutions obtained
by mirroring $\{0,1,\ldots,N\}\to \{N,N-1,\ldots,0\}$ are explored in
\S\,\ref{sec:dual}. Several semi-infinite Markov chains involving the little
$q$-Jacobi polynomials are derived in \S\,\ref{sec:lqJ} through dual
convolutions based on $q$-Hahn polynomials. The repeated discrete time birth
and death processes are presented in \S\,\ref{sec:repBD}.
The final section is for comments on the salient properties of the eigenvalues
of the $K$'s constructed in this project.
The proof of the triangularity {\bf Lemma} is provided as Appendix.

\section{Main Theorems}
\label{sec:mainT}

\subsection{Problem setting}
\label{sec:set}

We discuss stationary Markov chains on a one-dimensional finite integer
lattice $\cX$,
\begin{equation}
  \cX\eqdef\{0,1,\ldots,N\},
\end{equation}
with its points denoted by $x$, $y$, $z$, etc.
The main ingredient of the theory is the transition probability matrix $K(x,y)$
on $\cX$ which specifies the transition probability from $y$ to $x$
satisfying the basic conditions of probability and its conservation,
\begin{equation}
  K(x,y)\ge0,\quad\sum_{x\in\cX}K(x,y)=1.
  \label{basK}
\end{equation}
For given $K(x,y)$, much useful information can be extracted depending on the
specific needs of the applications.
The most basic ones would be the solutions of
the initial value problem and the determination of the $\ell$-step transition
probability:
\begin{align}
  &\bullet\text{\bf Initial value problem :}
  \ \ \text{one step time evolution}
  \ \cP(x;\ell+1)=\sum_{y\in\cX}K(x,y)\cP(y;\ell),\n %
  &\qquad\cP(x;0)\ge0,\ \sum_{x\in\cX}\cP(x;0)=1
  \Rightarrow\cP(x;\ell)=\sum_{y\in\cX}K^{\ell}(x,y)\cP(y;0),
  \label{inival}\\
  &\bullet\text{\bf $\ell$-step transition probability from $y$ to $x$ :}
  \ \cP(x,y;\ell),\n
  &\qquad\cP(x;0)=\delta_{x\,y}
  \Rightarrow\cP(x,y;\ell)=K^{\ell}(x,y).
  \label{ellpro}
\end{align}
They can be obtained based on the {\em complete set of eigenvalues and the
corresponding eigenvectors} of $K(x,y)$. It is well known that $K(x,y)$ has
always a maximal eigenvalue 1 and the range of spectrum
\begin{equation}
  -1\le\text{The moduli of the eigenvalues of }K(x,y)\le1,
\end{equation}
as a consequence of the non-negativity, {\em i.e.}\ Perron-Frobenius theorem,
and probability conservation \eqref{basK}.

\bigskip
Among many known strategies of procuring $K(x,y)$ with explicit forms of
eigensystems, one very promising plan is to construct $K(x,y)$ within the
framework of orthogonal polynomials on $\cX$. Simple examples are Birth
and Death (BD) processes, well known processes of nearest neighbour hopping.
All hypergeometric orthogonal polynomials of a discrete variable belonging to
Askey scheme provide exactly solvable continuous time BD processes
\cite{bdsol,os34} and a good part of it solves the discrete time versions
\cite{dtbd}, typical Markov chains. The normalised orthogonality measures
always provide the stationary probability distributions of the corresponding
BD and the polynomials the eigenvectors.
 
\subsection{Basic properties of $K(x,y)$}
\label{sec:base}

Here we develop a rather general method of building $K(x,y)$ by
`{\em convolutions}' of the normalised orthogonality measures. There are many
different types of convolutions available, but the main structure of the logic
is common, which we will describe by choosing one typical convolution.
In \S\,\ref{sec:exa} many explicit examples of $K(x,y)$ together with the
eigensystems, etc, constructed by various types of convolutions, will be
displayed together with multitudes of derivative forms obtained by the limiting
procedures of $N\to\infty$.

\bigskip
Let us introduce some notation. The normalised orthogonality measure of
orthogonal polynomials $\check{P}_n(x)\eqdef P_n(\eta(x))$
($\deg P_n(\eta)=n$) on $\cX$ is denoted by
$\pi(x,N,\bm{\lambda})$ ($x\in\cX$)
\begin{equation}
  \pi(x,N,\bm{\lambda})>0,\ \ \sum_{x\in\cX}\pi(x,N,\bm{\lambda})=1,
  \ \ \sum_{x\in\cX}\pi(x,N,\bm{\lambda})
  \check{P}_m(x,\bm{\lambda})\check{P}_n(x,\bm{\lambda})
  =0\ (m,n\in\cX,m\neq n),
  \label{pidef1}
\end{equation}
in which $\eta(x)$ is called the sinusoidal coordinate \cite{os12}.
For the explicit forms, see \eqref{triang2}. 
Here $\bm{\lambda}$ stands for the set of parameters other than the size
of the lattice $N$.
The $N$ dependence of the polynomials is suppressed for simplicity.
Throughout this paper we adopt the {\em universal normalisation} of the
polynomials,
\begin{equation}
  \check{P}_0(x,\bm{\lambda})=1\ \ (x\in\cX),\quad
  \check{P}_n(0,\bm{\lambda})=1\ \ (n\in\cX),
  \label{univnorm}
\end{equation}
and $\eta(0)=0$.
As mentioned above, later $\pi$ will denote the stationary probability
distribution of the Markov chain. The formulation presented in this subsection
would in principle apply
for any polynomials of a discrete variable \cite{nikiforov} in the Askey
family \cite{askey}--\cite{gasper}. But for the actual construction of $K(x,y)$
we deal with three kinds of orthogonal polynomials, the Krawtchouk (K), Hahn (H)
and $q$-Hahn ($q$H), and the initial is attached, {\em e.g.}\ 
$\pi_{\text{K}}(x,N,\bm{\lambda})$, $\check{P}_{\text{K}\,n}(x,\bm{\lambda})$
when formulas need specification of the polynomials.
In later subsections, more polynomials defined on a semi-infinite lattice
$\cX=\mathbb{Z}_{\ge0}$ will be discussed. They are the Charlier (C),
Meixner (M), and $q$-Meixner ($q$M), which are obtained
from (K), (H) and ($q$H) by limiting processes of some parameters in
$\bm{\lambda}$ and $N\to\infty$.

\bigskip
The rules of the game is to construct $K(x,y)$ by a convolution of two or more
$\pi(x)$ with certain choice of parameter dependence $\bm{\lambda}_1$ and
$\bm{\lambda}_2$, for example,
\begin{equation}
  \text{(\romannumeral1)}:\quad
  K(x,y)\eqdef\sum_{z=0}^{\min(x,y)}\pi(x-z,N-z,\bm{\lambda}_2)
  \pi(z,y,\bm{\lambda}_1)\quad(x,y\in\cX),
  \label{conv1}
\end{equation}
in such a way that $K(x,y)$ satisfies a symmetry condition with another $\pi$
with a parameter dependence $\bm{\lambda}$,
\begin{equation}
  K(x,y)\pi(y,N,\bm{\lambda})=K(y,x)\pi(x,N,\bm{\lambda})\quad(x,y\in\cX).
  \label{Kcond1}
\end{equation}
The explicit forms of five types of convolutions $K(x,y)$ are listed in
\S\,\ref{sec:convlist}.
It is easy to see that the above $K(x,y)$ \eqref{conv1} and four other types
listed in \S\,\ref{sec:convlist} satisfy the basic condition of probability
conservation $\sum_{x\in\cX}K(x,y)=1$ independently of the choices of
parameters $\bm{\lambda}_1$ and $\bm{\lambda}_2$ (and $\bm{\lambda}_3$).
That is, the probability conservation condition is satisfied irrespective of
the presence of $\pi(x,N,\bm{\lambda})$ for the symmetry condition.
Now $K$ is a {\em positive} matrix and its eigenvalues are greater than $-1$.
Let us note that the probability conservation $\sum_{x\in\cX}K(x,y)=1$ means
that $\check{\mathcal{V}}_0(x)\eqdef1$ is the left eigenvector of $K$ belonging
to the highest eigenvalue 1:
\begin{equation*}
  \sum_{y\in\cX}K(y,x)\check{\mathcal{V}}_0(y)=\check{\mathcal{V}}_0(x).
\end{equation*}
For obtaining the rest of left eigenvectors of $K$, the following {\bf Lemma} is
essential.\\
{\bf Lemma}
\,\,\textit{The transition matrices $K$'s generated
by five types of convolutions
\eqref{conv11}--\eqref{conv5} satisfy the triangularity condition}
\begin{align}
  \sum_{y\in\cX}K(y,x)\eta(y)^n&=\sum_{m=0}^na_{n\,m}\eta(x)^m\ \ (n\in\cX)
  \quad\bigl(a_{n\,m}=0\ \,\text{for}\ \,n<m\bigr),
  \label{triang1}\\
  \eta(x)&=\left\{
  \begin{array}{ll}
  x&:\text{(\romannumeral1)--(\romannumeral5)}\\
  q^{-x}-1&:\text{(\romannumeral1)},\,\text{(\romannumeral3)},
  \,\text{(\romannumeral4)}
  \end{array}\right..
  \label{triang2}
\end{align}
That is the r.h.s. of \eqref{triang1} is a polynomial in $\eta(x)$ of degree
at most $n$.
The vectors $\{\eta(x)^n\}$ ($n\in\cX$, $\eta(x)^0\eqdefrm1$) form a basis of
$\mathbb{R}^{N+1}$, because $\eta(x)$ is an increasing function.
Since the proof of {\bf Lemma} is of rather technical nature, it is consigned
to Appendix.
It should be noted that the triangularity of $K$ is the consequence of the
explicit forms of the convolutions \eqref{conv11}--\eqref{conv5} and it is
independent of the specific choice of the parameters $\bm{\lambda}_1$ and
$\bm{\lambda}_2$ (and $\bm{\lambda}_3$).

This {\bf Lemma} means that the left eigenvalue of $K$ is $\kappa(n)=a_{n\,n}$
and the corresponding left eigenvector is given by a certain degree $n$
polynomial in $\eta(x)$,
$\check{\mathcal{V}}_n(x)\eqdef\mathcal{V}_n(\eta(x))$
($\check{\mathcal{V}}_0(x)=\check{\mathcal{V}}_n(0)=1$),
\begin{equation}
  \sum_{y\in\cX}K(y,x)\check{\mathcal{V}}_n(y)
  =\kappa(n)\check{\mathcal{V}}_n(x)\quad(n\in\cX).
  \label{VK}
\end{equation}
By taking $y$ summation of the
symmetry condition \eqref{Kcond1}, we find that $\pi(x,N,\bm{\lambda})$ is the
eigenvector of $K(x,y)$ with the maximal eigenvalue 1,
\begin{equation*}
  \sum_{y\in\cX}K(x,y)\pi(y,N,\bm{\lambda})
  =\sum_{y\in\cX}K(y,x)\pi(x,N,\bm{\lambda})
  =\pi(x,N,\bm{\lambda}).
\end{equation*}
The set of all eigenvectors of $K$ is given by
\begin{equation}
  \sum_{y\in\cX}K(x,y)\pi(y,N,\bm{\lambda})\check{\mathcal{V}}_n(y)
  =\kappa(n)\pi(x,N,\bm{\lambda})\check{\mathcal{V}}_n(x)\quad(n\in\cX).
  \label{KV}
\end{equation}

Let us introduce the square root of the stationary distribution
$\pi(x,N,\bm{\lambda})$ and a diagonal matrix $\Phi$ on $\cX$,
\begin{equation}
  \hat{\phi}_0(x,\bm{\lambda})\eqdef\sqrt{\pi(x,N,\bm{\lambda})},\quad
  \Phi(x,x)\eqdef\hat{\phi}_0(x,\bm{\lambda}),\quad
  \Phi(x,y)\eqdef0\ \ (x\neq y).
  \label{Phidef}
\end{equation}
By a similarity transformation in terms of this $\Phi$, we define a matrix
$\mathcal{H}$
as follows:
\begin{equation}
  \mathcal{H}\eqdef\Phi^{-1}K\Phi,\quad
  \mathcal{H}(x,y)=\frac{\hat{\phi}_0(y,\bm{\lambda})}
  {\hat{\phi}_0(x,\bm{\lambda})}K(x,y).
  \label{KHrel}
\end{equation}
By dividing both sides of the symmetry condition \eqref{Kcond1} by
$\hat{\phi}_0(x,\bm{\lambda})\hat{\phi}_0(y,\bm{\lambda})$, we obtain
\begin{equation}
  \mathcal{H}(x,y)=\mathcal{H}(y,x)\quad(x,y\in\cX),
  \label{realsym}
\end{equation}
namely $\mathcal{H}$ is a real symmetric matrix. Hence it is diagonalizable
and its eigenvectors can be taken to be orthogonal with each other.
On the other hand, \eqref{VK} and \eqref{Kcond1} imply
\begin{equation}
  \sum_{y\in\cX}\mathcal{H}(x,y)\hat{\phi}_0(y,\bm{\lambda})
  \check{\mathcal{V}}_n(y)
  =\kappa(n)\hat{\phi}_0(x,\bm{\lambda})\check{\mathcal{V}}_n(x)
  \quad(n\in\cX),
\end{equation}
namely $\hat{\phi}_0(x,\bm{\lambda})\check{\mathcal{V}}_n(x)$'s are
eigenvectors of $\mathcal{H}$.
The polynomial $\check{\mathcal{V}}_n(x)$ of degree $n$ in $\eta(x)$,
being orthogonal with others with respect to the measure
$\hat{\phi}_0(x,\bm{\lambda})^2=\pi(x,N,\bm{\lambda})$, should be
$\check{P}_n(x,\bm{\lambda})$.
Therefore we obtain the following theorem.
\begin{theo}
\label{theo1}
The eigenvectors of the two matrices $\mathcal{H}=\Phi^{-1}K\Phi$ and $K$ are
described by the orthogonal polynomials
$\check{P}_n(x,\bm{\lambda})=P_n(\eta(x),\bm{\lambda})$ $(n\in\cX)$,
belonging to $\pi(x,N,\bm{\lambda})$,
\begin{align}
  \sum_{y\in\cX}\mathcal{H}(x,y)\hat{\phi}_0(y,\bm{\lambda})
  \check{P}_n(y,\bm{\lambda})
  &=\kappa(n)\hat{\phi}_0(x,\bm{\lambda})\check{P}_n(x,\bm{\lambda})\quad
  (n\in\cX),
  \label{Heig}\\
  \sum_{y\in\cX}K(x,y)\pi(y,N,\bm{\lambda})\check{P}_n(y,\bm{\lambda})
  &=\kappa(n){\pi}(x,N,\bm{\lambda})\check{P}_n(x,\bm{\lambda})\quad(n\in\cX).
  \label{Keig}
\end{align}
\end{theo}
Here we suppress the parameter dependence of the eigenvalues
$\{\kappa(n)\}$ for the simplicity of presentation.

The set of orthogonal vectors
$\{\hat{\phi}_0(x,\bm{\lambda})\check{P}_n(x,\bm{\lambda})\}$,
\begin{align}
  &\quad\sum_{x\in\cX}\hat{\phi}_0(x,\bm{\lambda})^2\check{P}_m(x,\bm{\lambda})
  \check{P}_n(x,\bm{\lambda})
  =\sum_{x\in\cX}\pi(x,N,\bm{\lambda})P_m\bigl(\eta(x),\bm{\lambda}\bigr)
  P_n\bigl(\eta(x),\bm{\lambda}\bigr)\n
  &=\sum_{x\in\cX}\pi(x,N,\bm{\lambda})\check{P}_m(x,\bm{\lambda})
  \check{P}_n(x,\bm{\lambda})
  =\frac{\delta_{m\,n}}{d_n^2}\quad(m,n\in\cX),\quad d_0=1,
  \label{Pnnorm}
\end{align}
form the complete set of eigenvectors of $\mathcal{H}$.
Now the scale of the polynomials $\{\check{P}_n(x,\bm{\lambda})\}$ is fixed
by the universal normalisation \eqref{univnorm}, the normalisation constant
$d_n^2$ is uniquely determined by the above formula.
It should be noted that, because of the context, the present definition of
$d_n^2$ corresponds to $d_n^2/d_0^2$ in our previous series of papers
\cite{bdsol,os34,os12}. The parameter dependence of $d_n^2$ is also suppressed.
Based on the universal normalisation condition of the polynomials
\eqref{univnorm}, we arrive at the {\em universal formula for the eigenvalues}
$\{\kappa(n)\}$ in terms of $K(x,y)$, $\pi(x,N,\bm{\lambda})$ and
$\{\check{P}_n(x,\bm{\lambda})\}$.
\begin{theo}
\label{theo2}
By setting $x=0$ in \eqref{Keig}, we obtain the universal expression of the
eigenvalues
\begin{equation}
  \kappa(n)=\sum_{y\in\cX}K(0,y)
  \frac{\pi(y,N,\bm{\lambda})}{\pi(0,N,\bm{\lambda})}\check{P}_n(y,\bm{\lambda})
  \ \ (n\in\cX),\quad\kappa(0)=1.
  \label{eigform}
\end{equation}
\end{theo}
Let us introduce the set of orthonormal eigenvectors of $\mathcal{H}$,
\begin{equation}
  \hat{\phi}_n(x,\bm{\lambda})\eqdef
  d_n\hat{\phi}_0(x,\bm{\lambda})\check{P}_n(x,\bm{\lambda}),\quad
  \sum_{x\in\cX}\hat{\phi}_m(x,\bm{\lambda})\hat{\phi}_n(x,\bm{\lambda})
  =\delta_{m\,n}\ \ (n,m\in\cX).
\end{equation}
\begin{theo}
\label{theo3}
The spectral representation of the real symmetric
matrix $\mathcal{H}$ provides that of $K$,
\begin{align}
  \mathcal{H}(x,y)&=\sum_{n\in\cX}\kappa(n)
  \hat{\phi}_n(x,\bm{\lambda})\hat{\phi}_n(y,\bm{\lambda}),\\
  K(x,y)&=\hat{\phi}_0(x,\bm{\lambda})\sum_{n\in\cX}\kappa(n)
  \hat{\phi}_n(x,\bm{\lambda})\hat{\phi}_n(y,\bm{\lambda})
  \cdot\hat{\phi}_0(y,\bm{\lambda})^{-1}\n
  &=\sum_{n\in\cX}\kappa(n)d_n^2\,\pi(x,N,\bm{\lambda})
  \check{P}_n(x,\bm{\lambda})\check{P}_n(y,\bm{\lambda}).
\end{align}
\end{theo}
\begin{theo}
\label{theo4}
The solution of the initial value problem of the Markov chain with the
transition rate $K(x,y)$ after $\ell$ steps is given by
\begin{equation}
  \cP(x;\ell)=\hat{\phi}_0(x,\bm{\lambda})\sum_{n\in\cX}c_n\kappa(n)^{\ell}
  \,\hat{\phi}_n(x,\bm{\lambda})
  =\pi(x,N,\bm{\lambda})\sum_{n\in\cX}c_nd_n\kappa(n)^{\ell}
  \check{P}_n(x,\bm{\lambda}),
  \label{Kinisol1}
\end{equation}
in which $\{c_n\}$ are determined by the expansion of the initial distribution
$\mathcal{P}(x;0)$,
\begin{align}
  &\quad\cP(x;0)=\hat{\phi}_0(x,\bm{\lambda})\sum_{n\in\cX}c_n
  \hat{\phi}_n(x,\bm{\lambda})\n
  &\Rightarrow
  c_n=\sum_{x\in\cX}\hat{\phi}_n(x,\bm{\lambda})
  \hat{\phi}_0(x,\bm{\lambda})^{-1}\cP(x;0)
  =d_n\sum_{x\in\cX}\check{P}_n(x,\bm{\lambda})\cP(x;0)\ \ (n\in\cX),
  \quad c_0=1.
  \label{cndef}
\end{align}
The $\ell$-step transition matrix from $y$ to $x$ is
\begin{align}
  \cP(x,y;\ell)=K^{\ell}(x,y)&=\hat{\phi}_0(x,\bm{\lambda})
  \sum_{n\in\cX}\kappa(n)^{\ell}\,\hat{\phi}_n(x,\bm{\lambda})
  \hat{\phi}_n(y,\bm{\lambda})\hat{\phi}_0(y,\bm{\lambda})^{-1}\n
  &=\pi(x,N,\bm{\lambda})\sum_{n\in\cX}d_n^2\kappa(n)^{\ell}
  \check{P}_n(x,\bm{\lambda})\check{P}_n(y,\bm{\lambda}).
  \label{Kellsol2}
\end{align}
Since $\kappa(0)=1$ and $-1<\kappa(n)<1$ $(n\neq0)$, the stationary distribution
is reached asymptotically,
\begin{equation}
  \lim_{\ell\to\infty}\cP(x;\ell)=\lim_{\ell\to\infty}\cP(x,y;\ell)
  =\pi(x,N,\bm{\lambda}).
\end{equation}
\end{theo}

It should be stressed that the results and theorems derived in this section are
valid for other choices of convolutions than \eqref{conv1} so long as the basic
condition \eqref{basK} and the symmetry condition \eqref{Kcond1} 
and Lemma are satisfied.

When a good $N\to\infty$ limit exists, it leads to a Markov chain on a
semi-infinite lattice $\cX=\mathbb{Z}_{\ge0}$, and the above theorems also
hold. That is, the symmetry condition \eqref{Kcond1}, the eigenvectors
\eqref{Keig} and the eigenvalue formula \eqref{eigform} are
\begin{align}
  &K(x,y)\pi(y,\bm{\lambda})=K(y,x)\pi(x,\bm{\lambda})\quad(x,y\in\cX),
  \label{Kcond1inf}\\
  &\sum_{y\in\cX}K(x,y)\pi(y,\bm{\lambda})\check{P}_n(y,\bm{\lambda})
  =\kappa(n){\pi}(x,\bm{\lambda})\check{P}_n(x,\bm{\lambda})\quad(n\in\cX),
  \label{Keiginf}\\
  &\kappa(n)=\sum_{y\in\cX}K(0,y)
  \frac{\pi(y,\bm{\lambda})}{\pi(0,\bm{\lambda})}\check{P}_n(y,\bm{\lambda})
  \ \ (n\in\cX),\quad\kappa(0)=1,
  \label{eigforminf}
\end{align}
where $\cX=\mathbb{Z}_{\ge0}$.

\section{Fundamental Data}
\label{sec:fund}
Here we present fundamental data for constructing and displaying the explicit
forms of various realisations of $K(x,y)$.
Starting with the list of convolutions in \S\,\ref{sec:convlist}, the basic data
of `orthogonality measures $=$ stationary distributions,' polynomials and the
normalisation constants $d_n^2$ etc are presented in \S\,\ref{sec:data}.

\subsection{List of convolutions}
\label{sec:convlist}

Here we list five forms of `convolutions' used for the construction of $K(x,y)$.
The list is not exhaustive at all.
A new and interesting convolution might be added in future.
\begin{align}
  \text{(\romannumeral1)}:&\ \ K(x,y)\eqdef\sum_{z=0}^{\min(x,y)}
  \!\pi(x-z,N-z,\bm{\lambda}_2)\pi(z,y,\bm{\lambda}_1),
  \label{conv11}\\
  \text{(\romannumeral2)}:&\ \ K(x,y)\eqdef\sum_{z=\max(0,x+y-N)}^{\min(x,y)}
  \!\!\!\!\!\!\!\!\!\pi(x-z,N-y,\bm{\lambda}_2)\pi(z,y,\bm{\lambda}_1),
  \label{conv2}\\
  \text{(\romannumeral3)}:&\ \ K(x,y)\eqdef\sum_{z=\max(x,y)}^N
  \!\!\!\!\pi(x,z,\bm{\lambda}_2)\pi(z-y,N-y,\bm{\lambda}_1),
  \label{conv3}\\
  \text{(\romannumeral4)}:&\ \ K(x,y)\eqdef\sum_{z_2=0}^{\min(x,y)}
  \!\!\pi(z_2,y,\bm{\lambda}_1)\!\!\!\sum_{z_1=\max(x,y)}^N
  \!\!\!\!\pi(x-z_2,z_1-z_2,\bm{\lambda}_3)\pi(z_1-y,N-y,\bm{\lambda}_2),
  \label{conv4}\\
  \text{(\romannumeral5)}:&\ \ K(x,y)\eqdef\sum_{z_2=0}^{\min(x,y)}
  \!\!\pi(z_2,y,\bm{\lambda}_1)\!\!\!\sum_{z_1=x+y-z_2}^N
  \!\!\!\!\pi(x-z_2,z_1-y,\bm{\lambda}_3)\pi(z_1-y,N-y,\bm{\lambda}_2).
  \label{conv5}
\end{align}
It is easy to convince oneself that the basic condition \eqref{basK} is
satisfied for each convolution.

It is obvious that these are very different from the standard forms of
convolutions, {\em e.g.}
\begin{equation*}
  (f\ast g)(x)=\sum_{z\in\cX}f(x-z)g(z),
\end{equation*}
since the formulas \eqref{conv11}--\eqref{conv5} must contain $x$ and $y$.
The above forms could be considered as deformations of convolutions containing
$x$ and $y$, like
\begin{equation*}
  (f\ast g)(x,y)=\sum_{z\in\cX}f(x-z)g(z-y).
\end{equation*}
Similar expressions appear in section \ref{sec:exa} during the reduction,
$N\to\infty$, processes.

The stationary probability measures $\pi(x,N,\bm{\lambda})$ presented in the
subsequent subsection have a remarkable property of `self-similarity' under
`ordinary' convolutions. This will be demonstrated in \S\,\ref{sec:selfsim}.

\subsection{Polynomials data}
\label{sec:data}

Here we provide the basic data of the participating polynomials.
Most are known facts collected for the consistency of notation, which is
standard.
For the explicit definitions of the basic quantities, {\em e.g.}\ $(a)_n$,
$(a\,;q)_n$, ${}_rF_s$ and ${}_r\phi_s$, consult \cite{askey,ismail}.
The data of the second family of orthogonal polynomials of the $q$-Meixner
\eqref{qMpim}--\eqref{qMort2} are not reported in standard references
\cite{nikiforov}--\cite{gasper}.
We believe some explicit expressions of the general formulas
\eqref{pieta1}--\eqref{pieta3}, {\em e.g.}\ \eqref{Ktrif1}, \eqref{Ktrif2}
etc.\ are new.
The parametrisation of some polynomials \cite{os12}, (H) and ($q$H), is
different from the conventional ones. There are many equivalent and different
looking expressions. We adopt the ones easy to grasp and simple to use.
Recall that $d_n>0$.

The measure $\pi(x,N,\bm{\lambda})$ is defined for $N\in\mathbb{Z}_{\ge0}$ and
$x\in\{0,1,\ldots,N\}$. For simplicity in presentation, we extend the domain
of definition to $x,N\in\mathbb{Z}$ by setting $\pi(x,N,\bm{\lambda})=0$
for otherwise.
Similarly, the domain of definition of $\pi(x,\bm{\lambda})$
($x\in\mathbb{Z}_{\ge0}$) is extended to $x\in\mathbb{Z}$ by setting
$\pi(x,\bm{\lambda})=0$ for otherwise.

\subsubsection{Krawtchouk (K)}

The polynomial depends on one positive parameter $\bm{\lambda}=p$ ($0<p<1$),
\begin{align}
  &\pi(x,N,p)=\binom{N}{x}p^x(1-p)^{N-x},\quad
  \binom{N}{x}=\frac{N!}{x!\,(N-x)!}\,,\quad
  d_n^2=\binom{N}{n}\Bigl(\frac{p}{1-p}\Bigr)^n,\\
  &\pi(N-x,N,p)=\pi(x,N,1-p),
  \label{Kpid}\\
  &s_1\eta(x)\pi(x,N,p)=-\pi(x-1,N-1,p),\quad\eta(x)=x,\quad
  s_1\eqdef-\frac{1}{pN},
  \label{Ktrif1}\\
  &\eta(z)\pi(z,x,p)=p\eta(x)\pi(z-1,x-1,p),
  \label{Ktrif2}\\
  &\check{P}_n(x,p)=P_n(x,p)
  ={}_2F_1\Bigl(\genfrac{}{}{0pt}{}{-n,\,-x}{-N}\Bigm|p^{-1}\Bigr),\quad
  P_n(x,p)=P_x(n,p),
  \label{Kp}\\
  &\check{P}_n(N-x,p)=(-1)^n(p^{-1}-1)^n\check{P}_n(x,1-p).
  \label{Kpd}
\end{align}
${P}_n(x,p)$ is a self-dual $P_n(x,p)=P_x(n,p)$ \cite{os12} degree $n$
polynomial in $x$ and $\pi$ is the binomial distribution.

\subsubsection{Charlier (C)}

This polynomial is defined on a semi-infinite integer lattice
$\cX=\mathbb{Z}_{\ge0}$ with $\bm{\lambda}=a$ ($a>0$),
\begin{align}
  &\pi(x,a)=\frac{a^xe^{-a}}{x!},\quad
  d_n^2=\frac{a^n}{n!},
  \label{Cpid}\\
  &s_1\eta(x)\pi(x,a)=-\pi(x-1,a),\quad\eta(x)=x,\quad s_1\eqdef-\frac{1}{a},
  \label{Ctrif1}\\
  &\check{P}_n(x,a)=P_n(x,a)
  ={}_2F_0\Bigl(\genfrac{}{}{0pt}{}{-n,\,-x}{-}\Bigm|-a^{-1}\Bigr),\quad
  P_n(x,a)=P_x(n,a).
  \label{Cp}
\end{align}
${P}_n(x,a)$ is a degree $n$ polynomial in $x$ and $\pi$ is the Poisson
distribution. It is self-dual $P_n(x,a)=P_x(n,a)$, too.
By the replacement $p\to pN^{-1}$ and the limit $N\to\infty$, the
Krawtchouk (K) goes to Charlier (C),
\begin{equation*}
  \check{P}_{\text{K}\,n}(x,p)\to\check{P}_{\text{C}\,n}(x,p),\quad
  \pi_{\text{K}}(x,N,p)\to\pi_{\text{C}}(x,p),\quad
  d_{\text{K}\,n}^2\to d_{\text{C}\,n}^2.
\end{equation*}

\subsubsection{Hahn (H)}

The polynomial depends on two positive parameters $\bm{\lambda}=(a,b)$
($a,b>0$),
\begin{align}
  &\pi(x,N,a,b)=\binom{N}{x}\frac{(a)_x\,(b)_{N-x}}{(a+b)_N},\quad
  d_n^2=\binom{N}{n}\frac{(a)_n\,(2n+a+b-1)(a+b)_N}{(b)_n\,(n+a+b-1)_{N+1}},\\
  &\pi(N-x,N,a,b)=\pi(x,N,b,a),
  \label{Hpi}\\
  &s_1\eta(x)\pi(x,N,a,b)=-\pi(x-1,N-1,a+1,b),\quad\eta(x)=x,\quad
  s_1\eqdef-\frac{a+b}{aN},
  \label{Htrif1}\\
  &\eta(z)\pi(z,x,a,b)=\frac{a}{a+b}\eta(x)\pi(z-1,x-1,a+1,b),
  \label{Htrif2}\\
  &\check{P}_n(x,a,b)=P_n(x,a,b)
  ={}_3F_2\Bigl(\genfrac{}{}{0pt}{}{-n,\,n+a+b-1,\,-x}{a,\,-N}\Bigm|1\Bigr),
  \label{Hp}\\
  &\check{P}_n(N-x,a,b)=(-1)^n\frac{(b)_n}{(a)_n}\check{P}_n(x,b,a).
  \label{Hpd}
\end{align}
${P}_n(x,a,b)$ is a degree $n$ polynomial in $x$ and $\pi$ is connected with
the hypergeometric distribution or the Polya distribution.

\subsubsection{Meixner (M)}

This polynomial is defined on a semi-infinite integer lattice
$\cX=\mathbb{Z}_{\ge0}$ with $\bm{\lambda}=(a,b)$ ($a>0$, $0<b<1$),
\begin{align}
  &\pi(x,a,b)=\frac{(a)_x\,b^x(1-b)^a}{x!},\quad
  d_n^2=\frac{(a)_n\,b^n}{n!},
  \label{Mpi}\\
  &s_1\eta(x)\pi(x,a,b)=-\pi(x-1,a+1,b),\quad\eta(x)=x,\quad
  s_1\eqdef-\frac{b^{-1}-1}{a},
  \label{Mtrif1}\\
  &\check{P}_n(x,a,b)=P_n(x,a,b)
  ={}_2F_1\Bigl(\genfrac{}{}{0pt}{}{-n,\,-x}{a}\Bigm|1-b^{-1}\Bigr),\quad
  P_n(x,a,b)=P_x(n,a,b).
  \label{Mp}
\end{align}
${P}_n(x,a,b)$ is a self-dual degree $n$ polynomial in $x$ and $\pi$ is
connected with the negative binomial distribution.
By the replacement $b\to N(1-b)b^{-1}$ and the limit $N\to\infty$, the Hahn (H)
goes to Meixner (M),
\begin{equation*}
  \check{P}_{\text{H}\,n}(x,a,b)\to\check{P}_{\text{M}\,n}(x,a,b),\quad
  \pi_{\text{H}}(x,N,a,b)\to\pi_{\text{M}}(x,a,b),\quad 
  d_{\text{H}\,n}^2\to d_{\text{M}\,n}^2.
\end{equation*}
By the replacement $b\to b/(a+b)$ and the limit $a\to\infty$, the Meixner
(M) goes to Charlier (C)
\begin{equation*}
  \check{P}_{\text{M}\,n}(x,a,b)\to\check{P}_{\text{C}\,n}(x,b),\quad
  \pi_{\text{M}}(x,a,b)\to\pi_{\text{C}}(x,b),\quad 
  d_{\text{M}\,n}^2\to d_{\text{C}\,n}^2.
\end{equation*}

\subsubsection{$q$-Hahn ($q$H)}

The three polynomials, the $q$-Hahn ($q$H), $q$-Meixner ($q$M) and $q$-Charlier
($q$C), to be discussed hereafter, depend on $q$, $0<q<1$ on top of the other
parameters. The $q$ dependence of $\pi$ and $\check{P}_n$ is suppressed.
The limiting processes of these $q$-polynomials to non $q$-polynomials will not
be discussed here.
It should be stressed that these three polynomials $\check{P}_n(x)$ are
degree $n$ polynomials in $q^{-x}-1$, not in $x$.
The $q$-Hahn is defined on a finite integer lattice with two positive parameters
$\bm{\lambda}=(a,b)$ ($0<a<1$, $b<1$),
\begin{align}
  &\pi(x,N,a,b)=\genfrac{[}{]}{0pt}{}{\,N\,}{x}
  \frac{(a\,;q)_x\,(b\,;q)_{N-x}a^{N-x}}{(ab\,;q)_N},\quad
  \genfrac{[}{]}{0pt}{}{\,N\,}{x}\eqdef\frac{(q\,;q)_N}
  {(q\,;q)_x\,(q\,;q)_{N-x}},
  \label{qHpi}\\
  &d_n^2=\genfrac{[}{]}{0pt}{}{\,N\,}{n}
  \frac{(a,abq^{-1}\,;q)_n}{(abq^N,b\,;q)_n\,a^n}\frac{1-abq^{2n-1}}{1-abq^{-1}},
  \label{qHd}\\
  &s_1\eta(x)\pi(x,N,a,b)=-\pi(x-1,N-1,aq,b),\quad
  s_1\eqdef-\frac{1-ab}{(1-a)(q^{-N}-1)},
  \label{qHtrif1}\\
  &\eta(z)\pi(z,x,a,b)=\frac{1-a}{1-ab}\eta(x)\pi(z-1,x-1,aq,b),\quad
  \eta(x)=q^{-x}-1,
  \label{qHtrif2}\\
  &\check{P}_n(x,a,b)=P_n\bigl(\eta(x),a,b\bigr)
  ={}_3\phi_2\Bigl(\genfrac{}{}{0pt}{}{q^{-n},\,abq^{n-1},\,q^{-x}}
  {a,\,q^{-N}}\Bigm|q\,;q\Bigr),
  \label{qHp}\\
  &\pi(N-x,N,a,b)=\frac{(ab)^x}{b^N}\pi(x,N,b,a),\quad
  \check{P}_n(N-x,N,a,b)\not\propto\check{P}_n(x,N,b,a).
  \label{qHrn}
\end{align}

\subsubsection{$q$-Meixner ($q$M)}

This is a polynomial in $\eta(x)=q^{-x}-1$ defined on a semi-infinite integer
lattice $\cX=\mathbb{Z}_{\ge0}$ with $\bm{\lambda}=(b,c)$ ($0<b<q^{-1}$, $c>0$),
\begin{align}
  &\pi(x,b,c)=\frac{(bq\,;q)_x}{(q,-bcq\,;q)_x}\,c^xq^{\binom{x}{2}}
  \frac{(-bcq\,;q)_{\infty}}{(-c\,;q)_{\infty}},\quad
  d_n^2=\frac{q^n(bq\,;q)_n}{(q,-c^{-1}q\,;q)_n},
  \label{qMpi}\\
  &s_1\eta(x)\pi(x,b,c)=-\pi(x-1,bq,c),\quad\eta(x)=q^{-x}-1,\quad
  s_1\eqdef-\frac{q}{c(1-bq)},
  \label{qMtrif11}\\
  &\check{P}_n(x,b,c)=P_n\bigl(\eta(x),b,c\bigr)
  ={}_2\phi_1\Bigl(
  \genfrac{}{}{0pt}{}{q^{-n},\,q^{-x}}{bq}\Bigm|q\,;-c^{-1}q^{n+1}\Bigr),\\
  &\hat{\phi}_0(x,b,c)=\sqrt{\pi(x,b,c)},\quad
  \hat{\phi}_n(x,b,c)=d_n\hat{\phi}_0(x,b,c)\check{P}_n(x,b,c),\n
  &\sum_{x\in\cX}\hat{\phi}_n(x,b,c)\hat{\phi}_m(x,b,c)=\delta_{n\,m}
  \ \ (n,m\in\cX).
  \label{qMort}
\end{align}
The completeness relation is not satisfied
\begin{equation*}
  \sum_{n\in\cX}\hat{\phi}_n(x,b,c)\hat{\phi}_n(y,b,c)\neq\delta_{x\,y}
  \ \ (x,y\in\cX),
\end{equation*}
by these polynomials \cite{atakishi1} as seen clearly by (3.15) of \cite{os34}.
Another set of orthogonal polynomials obtained from the original set by the
parameter change (involution)
\begin{equation*}
  (b,c)\to(-bc,c^{-1}),
\end{equation*}
is needed for the completeness,
\begin{align}
  &\pi^{(-)}(x,b,c)=\frac{(-bcq\,;q)_x}{(q,bq\,;q)_x}\,c^{-x}q^{\binom{x}{2}}
  \frac{(bq\,;q)_{\infty}}{(-c^{-1}\,;q)_{\infty}},
  \ \ d_n^{(-)\,2}=\frac{q^n(-bcq\,;q)_n}{(q,-cq\,;q)_n}\ (d_n^{(-)}>0),
  \label{qMpim}\\
  &s_1\eta(x)\pi^{(-)}(x,b,c)=-\pi^{(-)}(x-1,bq,c),\quad
  s_1\eqdef-\frac{cq}{1+bcq},
  \label{qMtrif21}\\
  &\check{P}^{(-)}_n(x,b,c)=P^{(-)}_n\bigl(\eta(x),b,c\bigr)
  ={}_2\phi_1\Bigl(\genfrac{}{}{0pt}{}
  {q^{-n},q^{-x}}{-bcq}\!\!\Bigm|\!q\,;-cq^{n+1}\Bigr),
  \label{qM:Pnm}\\
  &\hat{\phi}^{(-)}_0(x,b,c)\eqdef(-1)^x\sqrt{\pi^{(-)}(x,b,c)},\quad
  \hat{\phi}^{(-)}_n(x,b,c)\eqdef d_n^{(-)}\hat{\phi}^{(-)}_0(x,b,c)
  \check{P}^{(-)}_n(x,b,c),\\
  &\sum_{x\in\cX}\hat{\phi}_n(x,b,c)\hat{\phi}^{(-)}_m(x,b,c)=0,\quad
  \sum_{x\in\cX}\hat{\phi}^{(-)}_n(x,b,c)\hat{\phi}^{(-)}_m(x,b,c)=\delta_{n\,m}
  \ \ (n,m\in\cX).
 \label{qMort2}
\end{align}
$q$M \eqref{qMpi}--\eqref{qMort2} is obtained from $q$H by the replacement
$a\to bq$, $b\to-b^{-1}c^{-1}q^{-N}$ and the limit $N\to\infty$.
The $q$-Charlier with $\bm{\lambda}=c$ ($c>0$) is obtained from $q$-Meixner by
setting $b=0$.

\subsection{Convolutional self-similarity of stationary distributions}
\label{sec:selfsim}

It is well known that the Gaussian distribution
\begin{equation*}
  \pi_{\text{G}}(x,\sigma)\eqdef\frac1{\sqrt{2\pi}\,\sigma}\,
  e^{-\tfrac{x^2}{2\sigma^2}}\quad(\sigma>0),
\end{equation*}
keeps its form under the standard convolution
\begin{equation*}
  \int_{-\infty}^{\infty}\pi_{\text{G}}(x-z,\sigma)\pi_{\text{G}}(z,\tau)dz
  =\pi_{\text{G}}(x,\sqrt{\sigma^2+\tau^2}\,).
\end{equation*}
Here we show that the $\pi$'s listed in the previous subsection also keep their
forms under several forms of convolutions. This means repeating these
convolutions as a part of constructing $K(x,y)$ would be redundant, as using
two $\pi$'s is the same as one $\pi$.
The following formulas are verified easily by straightforward calculation.

\paragraph{Charlier (C)}

\begin{align}
  \sum_{z=0}^x\pi(x-z,a_2)\pi(z,a_1)&=\pi(x,a_1+a_2),
  \label{CC1}\\
  \sum_{z=y}^x\pi(x-z,a_2)\pi(z-y,a_1)&=\pi(x-y,a_1+a_2).
  \label{CC2}
\end{align}
These are classical results obtained by the binomial theorem.

\paragraph{Meixner (M)}

\begin{align}
  \sum_{z=0}^x\pi(x-z,a_2,b)\pi(z,a_1,b)&=\pi(x,a_1+a_2,b),
  \label{MM1}\\
  \sum_{z=y}^x\pi(x-z,a_2,b)\pi(z-y,a_1,b)&=\pi(x-y,a_1+a_2,b).
  \label{MM2}
\end{align}
These are obtained by the sum formula of $\pi_{\text{H}}$
\begin{equation}
  \sum_{x\in\cX}\pi_{\text{H}}(x,N,a,b)=1\Longleftrightarrow
  \sum_{k=0}^n\binom{n}{k}(a)_k(b)_{n-k}=(a+b)_n\ \ (n\in\mathbb{Z}_{\ge0}),
  \label{bt2}
\end{equation}
which is derived from the formula \cite{kls}(1.5.4),
\begin{equation}
  {}_2F_1\Bigl(\genfrac{}{}{0pt}{}{-n,\,b}{c}\Bigm|1\Bigr)
  =\frac{(c-b)_n}{(c)_n}\ \ (n\in\mathbb{Z}_{\ge0}).
  \label{KLS(1.5.4)}
\end{equation}
These four formulas (C) and (M) are symmetric in $\bm{\lambda}_1$ and
$\bm{\lambda}_2$.

\paragraph{Krawtchouk (K)}

\begin{align}
  \sum_{z=0}^x\pi(x-z,N-z,a_2)\pi(z,N,a_1)
  &=\pi\bigl(x,N,1-(1-a_1)(1-a_2)\bigr),
  \label{thconK}\\
  \sum_{z=x}^y\pi(x,z,a_2)\pi(z,y,a_1)&=\pi(x,y,a_1a_2),
  \label{Kpipi(i)}\\
  \sum_{z=y}^x\pi(x-z,N-z,a_2)\pi(z-y,N-y,a_1)
  &=\pi\bigl(x-y,N-y,1-(1-a_1)(1-a_2)\bigr).
  \label{Kpipi(ii)}
\end{align}
In all formulas the binomial theorem is used.
The result \eqref{thconK} was reported in \cite{hoa-rah83}p115, (2.3) together
with the $n$-fold repetition in (2.6).
All three formulas work when $\bm{\lambda}_1=a_1$ and $\bm{\lambda}_2=a_2$ are
interchanged.

\paragraph{Hahn (H)}

\begin{align}
  \sum_{z=0}^x\pi(x-z,N-z,a_1,b_1)\pi(z,N,a_2,a_1+b_1)&=\pi(x,N,a_1+a_2,b_1),\\
  \sum_{z=x}^y\pi(x,z,a_1,b_1)\pi(z,y,a_1+b_1,b_2)&=\pi(x,y,a_1,b_1+b_2),
  \label{hconv2}\\
  \sum_{z=y}^x\pi(x-z,N-z,a_1,b_1)\pi(z-y,N-y,a_2,a_1+b_1)
  &=\pi(x-y,N-y,a_1+a_2,b_1).
  \label{hconv3}
\end{align}
These formulas are obtained by the sum formula of $\pi_{\text{H}}$ \eqref{bt2}.
It is very interesting to note that, like (C), (M) and (K), \eqref{hconv2} and
\eqref{hconv3} work when $\bm{\lambda}_1$ and $\bm{\lambda}_2$ are interchanged,
\begin{align}
  \sum_{z=x}^y\pi(x,z,a_1+b_1,b_2)\pi(z,y,a_1,b_1)&=\pi(x,y,a_1,b_1+b_2),
  \label{hconv4}\\
  \sum_{z=y}^x\pi(x-z,N-z,a_2,a_1+b_1)\pi(z-y,N-y,a_1,b_1)
  &=\pi(x-y,N-y,a_1+a_2,b_1).
  \label{hconv5}
\end{align}
These are obtained by the sum formula ($n,m\in\mathbb{Z}_{\ge0}$),
\begin{equation}
  \sum_{k=0}^n\binom{n}{k}
  (b_1)_{n-k}(b_2)_k\frac{(a)_{m+k}}{(a+b_1+b_2)_{m+k}}
  =\frac{(a)_m(b_1+b_2)_n}{(a+b_1+b_2)_{m+n}}
  \frac{(a+b_1)_{m+n}}{(a+b_1)_{m}},
  \label{bt3}
\end{equation}
which is derived from the Pfaff-Saalsch\"utz formula \cite{askey}(2.2.8),
\begin{equation}
  {}_3F_2\Bigl(
  \genfrac{}{}{0pt}{}{-n,\,a,\,b}{c,\,1+a+b-c-n}\Bigm|1\Bigr)
  =\frac{(c-a,c-b)_n}{(c,c-a-b)_n}\ \ (n\in\mathbb{Z}_{\ge0}).
  \label{KLS(1.5.5)}
\end{equation}

\paragraph{$q$-Hahn ($q$H)}

\begin{align}
  \sum_{z=0}^x\pi(x-z,N-z,a_1,b_1)\pi(z,N,a_2,a_1b_1)
  &=\pi(x,N,a_1a_2,b_1),\\
  \sum_{z=x}^y\pi(x,z,a_1,b_1)\pi(z,y,a_1b_1,b_2)&=\pi(x,y,a_1,b_1b_2),
  \label{qhconv2}\\
  \sum_{z=y}^x\pi(x-z,N-z,a_1,b_1)\pi(z-y,N-y,a_2,a_1b_1)
  &=\pi(x-y,N-y,a_1a_2,b_1).
\label{qhconv3}
\end{align}
These are obtained by the sum formula for $\pi_{\text{$q$H}}$
\begin{equation}
  \sum_{x\in\cX}\pi_{\text{$q$H}}(x,N,a,b)=1\Longleftrightarrow
  \sum_{k=0}^n\genfrac{[}{]}{0pt}{}{\,n\,}{k}(a\,;q)_k(b\,;q)_{n-k}a^{n-k}
  =(ab\,;q)_n\ \ (n\in\mathbb{Z}_{\ge0}),
  \label{bt4}
\end{equation}
which is derived from the formula \cite{kls}(1.11.4),
\begin{equation}
  {}_2\phi_1\Bigl(\genfrac{}{}{0pt}{}{q^{-n},\,b}{c}\Bigm|
  q\,;\frac{cq^n}{b}\Bigr)
  =\frac{(b^{-1}c\,;q)_n}{(c\,;q)_n}\ \ (n\in\mathbb{Z}_{\ge0}).
  \label{KLS(1.11.4)}
\end{equation}
Similarly to the Hahn cases, \eqref{qhconv2} and \eqref{qhconv3} work when
$\bm{\lambda}_1$ and $\bm{\lambda}_2$ are interchanged,
\begin{align}
  \sum_{z=x}^y\pi(x,z,a_1b_1,b_2)\pi(z,y,a_1,b_1)&=\pi(x,y,a_1,b_1b_2),
  \label{qhconv4}\\
  \sum_{z=y}^x\pi(x-z,N-z,a_2,a_1b_1)\pi(z-y,N-y,a_1,b_1)
  &=\pi(x-y,N-y,a_1a_2,b_1).
  \label{qhconv5}
\end{align}
These are obtained by the sum formula ($n,m\in\mathbb{Z}_{\geq 0}$)
\begin{equation}
  \sum_{k=0}^n\genfrac{[}{]}{0pt}{}{\,n\,}{k}(b_1\,;q)_{n-k}(b_2\,;q)_k
  \frac{b_1^k\,(a\,;q)_{m+k}}{(ab_1b_2\,;q)_{m+k}}
  =\frac{(a\,;q)_m(b_1b_2\,;q)_n}{(ab_1b_2\,;q)_{m+n}}
  \frac{(ab_1\,;q)_{m+n}}{(ab_1\,;q)_m},
  \label{bt5}
\end{equation}
which is derived from the $q$-Pfaff-Saalsch\"utz formula \cite{kls}(1.11.9),
\begin{equation}
  {}_3\phi_2\Bigl(
  \genfrac{}{}{0pt}{}{q^{-n},\,a,\,b}{c,\,abc^{-1}q^{1-n}}\Bigm|q\,;q\Bigr)
  =\frac{(a^{-1}c,b^{-1}c;q)_n}{(c,a^{-1}b^{-1}c;q)_n}
  \ \ (n\in\mathbb{Z}_{\ge0}).
  \label{KLS(1.11.9)}
\end{equation}

\section{Many Examples of $K(x,y)$}
\label{sec:exa}

In this section we present various examples of $K(x,y)$ constructed by the five
types of convolutions listed in \S\,\ref{sec:convlist} applied to the
polynomials Krawtchouk (K), Hahn (H) and $q$-Hahn ($q$H).
For each polynomial, the limiting forms obtained by $N\to\infty$ are displayed,
Charlier (C), Meixner (M) and $q$-Meixner ($q$M).
We do strongly hope that these examples would enrich many related disciplines.

The symmetry condition \eqref{Kcond1} or \eqref{Kcond1inf} is easily verified
without evaluating the sum(s) in $K(x,y)$ for all examples in this section
except for those in \S\,\ref{sec:(i)(iii)gen}.

\subsection{Type (\romannumeral1) convolution}
\label{convtype1}

This convolution 
\begin{equation*}
  \text{(\romannumeral1)}:\ \ K(x,y)\eqdef\sum_{z=0}^{\min(x,y)}
  \!\pi(x-z,N-z,\bm{\lambda}_2)\pi(z,y,\bm{\lambda}_1),
\end{equation*}
has been applied to (K) and (H) in many papers
\cite{hoa-rah83}--\cite{albert} in connection with ``cumulative Bernoulli
trials.'' It has a factorised form,
$K(x,y)=\sum_{z=0}^NA(x,z)B(z,y)$, the non-vanishing elements
of $A(x,z)$ and $B(z,y)$ being $\pi(x-z,N-z,\bm{\lambda}_2)$ ($x\ge z$) and
$\pi(z,y,\bm{\lambda}_1)$ ($z\le y$), respectively.
Namely, $A$ is a lower triangular matrix and $B$ is upper triangular.
The determinant of $K$ is easily obtained as
\begin{equation*}
  \prod_{n\in\cX}\kappa(n)=\det K
  =\prod_{x\in\cX}A(x,x)\cdot\prod_{x\in\cX}B(x,x)
  =\prod_{n\in\cX}\pi(0,n,\bm{\lambda}_2)\pi(n,n,\bm{\lambda}_1),
\end{equation*}
from which eigenvalues $\kappa(n)$ are obtained, if it is known that the
eigenvalues are independent of $N$
($\Rightarrow$ $\kappa(n)=\pi(0,n,\bm{\lambda}_2)\pi(n,n,\bm{\lambda}_1)$).
In fact, as we will see shortly, the eigenvalues are $N$ independent for all
the examples in this section.
Therefore for type (\romannumeral1) and (\romannumeral3) convolutions, the
determinant formulas give eigenvalues.
Moreover, for $x=0$ only one term $z=0$ contributes to $K(0,y)$. This greatly
simplifies the eigenvalue formula \eqref{eigform}
\begin{equation}
  \kappa(n)=\sum_{y\in\cX}\pi(0,N,\bm{\lambda}_2)\pi(0,y,\bm{\lambda}_1)
  \frac{\pi(y,N,\bm{\lambda})}{\pi(0,N,\bm{\lambda})}\check{P}_n(y,\bm{\lambda})
  \ \ (n\in\cX).
  \label{eigform1}
\end{equation}

\subsubsection{Krawtchouk (K)}

By taking $\bm{\lambda}_1=a$ and $\bm{\lambda}_2=b$, the matrix $K(x,y)$ is
\begin{equation}
 K(x,y)=\!\sum_{z=0}^{\min(x,y)}\pi(x-z,N-z,b)\pi(z,y,a)\quad
  \Bigl(\Rightarrow\det K=\prod_{n\in\cX}a^n(1-b)^n\Bigr).
  \label{KK1}
\end{equation}
For the following $\bm{\lambda}$,
\begin{equation}
  \bm{\lambda}=p\eqdef\frac{b}{1-a+ab},
  \label{K1la}
\end{equation}
the symmetry condition \eqref{Kcond1} is satisfied and
{\bf Theorem\,\ref{theo1}} gives \eqref{Keig},
\begin{equation*}
  \sum_{y\in\cX}K(x,y)\pi(y,N,p)\check{P}_n(y,p)
  =\kappa(n)\pi(x,N,p)\check{P}_n(x,p).
\end{equation*}
By writing down the eigenvalue formula \eqref{eigform1}, we have
\begin{equation}
  \kappa(n)=\sum_{y\in\cX}\pi(y,N,b)\check{P}_n(y,p)=a^n(1-b)^n,
  \label{Kkapi}
\end{equation}
by using a generating function of Krawtchouk $P_n(x)$ \cite{kls}(9.11.11),
\begin{equation}
  \sum_{n=0}^N\binom{N}{n}P_n(x,p)t^n
  =\Bigl(1-\frac{1-p}{p}t\Bigr)^x(1+t)^{N-x},
  \label{KLS(9.11.11)}
\end{equation}
together with the self-duality of (K) $P_n(x,p)=P_x(n,p)$.
In a summation formula like \eqref{Kkapi} $y$ summation is easily performed as 
the only $y$ dependence of $\check{P}_n(y,\bm{\lambda})$ in this paper is due
to $(-y)_k$
(or $(q^{-y}\,;q)_k$ for $q$H, $q$M and $q$C) ($k=0,1,\ldots,n$), which gives
\begin{equation}
  \sum_{y\in\cX}\pi(y,N,b)(-y)_k=(-N)_kb^k.
  \label{Nind}
\end{equation}
This cancels the $N$ dependence in the hypergeometric summation of
$\check{P}_n(y,p)$. This is the general mechanism leading to the
$N$-independence of the eigenvalues for (K), (H) and ($q$H).
By direct calculation using \eqref{Nind}, we obtain
\begin{equation}
  \kappa(n)=\sum_{k=0}^n\frac{(-n)_k\,p^{-k}}{(-N)_k\,k!}(-N)_kb^k
  =\sum_{k=0}^n(-n)_k\frac{(bp^{-1})^k}{k!}
  ={}_1F_0\Bigl(\genfrac{}{}{0pt}{}{-n}{-}\Bigm|bp^{-1}\Bigr).
\end{equation}
As we will show in the following, {\em all the eigenvalues $\kappa(n)$ of $K$'s
constructed in this section are expressed by a terminating
($q$-)hypergeometric series}, except for those corresponding to the extra
eigenvectors of ($q$M).

\paragraph{Krawtchouk\,$\to$\,Charlier}

This is achieved by $b\to bN^{-1}$, $N\to\infty$,
\begin{align}
  &\check{P}_n(x,p)\to\check{P}_{\text{C}\,n}(x,p'),\quad
  p'\eqdef\frac{b}{1-a},\n
  &\pi(x,N,p)\to\pi_{\text{C}}(x,p'),\quad
  \kappa(n)\to\kappa_{\text{C}}(n)=a^n,\n
  &K(x,y)\to K_{\text{C}}(x,y,a,b)
  =\sum_{z=0}^{\min(x,y)}\pi_{\text{C}}(x-z,b)\pi(z,y,a),
  \label{KCconv1}
\end{align}
and the relations \eqref{Kcond1}, \eqref{Keig} and \eqref{eigform} of (K)
reduce to those of (C).
The resulting $K_{\text{C}}$ is a standard convolution of the $\pi$'s of (C)
and (K). This was discussed in \cite{hoa-rah83}\S\,4.
Based on $K_{\text{C}}$ \eqref{KCconv1}, let us rederive these results.
The symmetry condition \eqref{Kcond1inf} is satisfied for $\bm{\lambda}=p'$ and
{\bf Theorem\,\ref{theo1}} gives \eqref{Keiginf}.
The direct calculation of the eigenvalue formula \eqref{eigforminf} gives
\begin{equation}
  \kappa_{\text{C}}(n)
  =\sum_{y=0}^{\infty}\pi_{\text{C}}(y,b)\check{P}_{\text{C}\,n}(y,p')=a^n= 
  {}_1F_0\Bigl(\genfrac{}{}{0pt}{}{-n}{-}\Bigm|bp^{\prime\,-1}\Bigr).
\end{equation}
It can also be obtained by using a generating function of (C) \cite{kls}(9.14.11)
together with the self-duality of $(C)$.

\subsubsection{Hahn (H)}

By taking $\bm{\lambda}_1=(a,b)$ and $\bm{\lambda}_2=(b,c)$, the matrix
$K(x,y)$ is
\begin{equation}
  K(x,y)=\!\sum_{z=0}^{\min(x,y)}\pi(x-z,N-z,b,c)\pi(z,y,a,b)
  \ \ \Bigl(\Rightarrow\det K=\prod_{n\in\cX}\frac{(a)_n(c)_n}{(a+b)_n(b+c)_n}
  \Bigr).
  \label{h12K}
\end{equation}
For the following $\bm{\lambda}$,
\begin{equation}
  \bm{\lambda}=(a+b,c),
  \label{H1la}
\end{equation}
the symmetry condition \eqref{Kcond1} is satisfied and
{\bf Theorem\,\ref{theo1}} gives \eqref{Keig},
\begin{equation*}
  \sum_{y\in\cX}K(x,y)\pi(y,N,a+b,c)\check{P}_n(y,a+b,c)
  =\kappa(n)\pi(x,N,a+b,c)\check{P}_n(x,a+b,c).
\end{equation*}
By evaluating the eigenvalue formula \eqref{eigform1}, we obtain a balanced
${}_3F_2$,
\begin{align}
  \kappa(n)&=\sum_{y\in\cX}\pi(y,N,b,c)\check{P}_n(y,a+b,c)\n
  &={}_3F_2\Bigl(
  \genfrac{}{}{0pt}{}{-n,\,n+a+b+c-1,\,b}{a+b,\,b+c}\Bigm|1\Bigr)
  =\frac{(a)_n(c)_n}{(a+b)_n(b+c)_n}.
  \label{1HF}
\end{align}
The last equality is due to the Pfaff-Saalsch\"utz formula \eqref{KLS(1.5.5)}.
This provides another sum formula for the dual Hahn polynomial
$\check{Q}_n(x,a+b,c)$ \cite{os12},
\begin{align}
  &\sum_{n=0}^N\binom{N}{n}(b)_n(c)_{N-n}\check{Q}_n(x,a+b,c)
  =\frac{(b+c)_N(a)_x(c)_x}{(a+b)_x(b+c)_x},\\
  &\check{Q}_n(x,a+b,c)\eqdef{}_3F_2\Bigl(
  \genfrac{}{}{0pt}{}{-n,\,x+a+b+c-1,\,-x}{a+b,\,-N}\Bigm|1\Bigr).
\end{align}
By rewriting $(a)_x=\Gamma(a+x)/\Gamma(a)$, this sum formula is valid for
$\forall x\in\mathbb{C}$.

\paragraph{Hahn\,$\to$\,Meixner}

This is achieved by fixing $a$ and $b$ with $c\to N(1-c)c^{-1}$
($\Rightarrow 0<c<1$), $N\to\infty$,
\begin{align}
  &\check{P}_n(x,a+b,c)\to\check{P}_{\text{M}\,n}(x,a+b,c),\n
  &\pi(x,N,a+b,c)\to\pi_{\text{M}}(x,a+b,c),\quad
  \kappa(n)\to\kappa_{\text{M}}(n)=\frac{(a)_n}{(a+b)_n},\n
  &K(x,y)\to K_{\text{M}}(x,y,a,b,c)=\sum_{z=0}^{\min(x,y)}
  \pi_{\text{M}}(x-z,b,c)\pi(z,y,a,b),
  \label{HM1}
\end{align}
and the relations \eqref{Kcond1}, \eqref{Keig} and \eqref{eigform} of (H)
reduce to those of (M).
This $K_{\text{M}}$ is a standard convolution of $\pi$'s of (M) and (H).
Based on $K_{\text{M}}$ \eqref{HM1}, let us rederive these results.
The symmetry condition \eqref{Kcond1inf} is satisfied for
$\bm{\lambda}=(a+b,c)$ and
{\bf Theorem\,\ref{theo1}} gives \eqref{Keiginf}.
The direct calculation of the eigenvalue formula \eqref{eigforminf} gives
\begin{equation}
  \kappa_{\text{M}}(n)=\sum_{y=0}^{\infty}\pi_{\text{M}}(y,b,c)
  \check{P}_{\text{M}\,n}(y,a+b,c)=\frac{(a)_n}{(a+b)_n}
  ={}_2F_1\Bigl(
  \genfrac{}{}{0pt}{}{-n,\,b}{a+b}\Bigm|1\Bigr).
\end{equation}
It can be obtained by using a generating function \cite{kls}(9.10.13)
together with the self-duality of (M).
Of course, ${}_2F_1$ form is also obtained from the balanced
${}_3F_2$ form \eqref{1HF} by the above limit.

\paragraph{Hahn\,$\to$\,Meixner\,$\to$\,Charlier}

This is achieved by $a\to aN$, $b\to bN$, $c\to\frac{c}{c+N}$, $N\to\infty$,
\begin{align}
  &P_{\text{M}\,n}\bigl(x,(a+b)N,\tfrac{c}{c+N}\bigr)\to
  P_{\text{C}\,n}\bigl(x,(a+b)c\bigr),\n
  &\pi_{\text{M}}\bigl(x,(a+b)N,\tfrac{c}{c+N}\bigr)
  \to\pi_{\text{C}}\bigl(x,(a+b)c\bigr),\n
  &\kappa_{\text{M}}(n)=\frac{(aN)_n}{\bigl((a+b)N\bigr)_n}\to
  \kappa_{\text{C}}(n)=\Bigl(\frac{a}{a+b}\Bigr)^n,\n
  &\pi_{\text{M}}\bigl(x,bN,\tfrac{c}{c+N}\bigr)\to\pi_{\text{C}}(x,bc),\quad
  \pi(z,y,aN,bN)\to\pi_{\text{K}}\bigl(z,y,\tfrac{a}{a+b}\bigr),\n
  &K_{\text{M}}(x,y,a,b,c)\to K_{\text{C}}(x,y,a,b,c)
  =\sum_{z=0}^{\min(x,y)}\pi_{\text{C}}(x-z,bc)
  \pi_{\text{K}}\bigl(z,y,\tfrac{a}{a+b}\bigr),
  \label{KCK}
\end{align}
and the relations \eqref{Kcond1inf}, \eqref{Keiginf} and \eqref{eigforminf}
of (M) reduce to those of (C).
This $K_{\text{C}}$ agrees with \eqref{KCconv1} with the replacement
$(a,b)\to(\frac{a}{a+b},bc)$.

\subsubsection{$q$-Hahn ($q$H)}
\label{1qH}

We believe the explicit examples of Markov chains described by the $q$-Hahn
polynomial and its reduction, $q$-Meixner are new.

By taking $\bm{\lambda}_1=(a,b)$ and $\bm{\lambda}_2=(b,c)$, the matrix
$K(x,y)$ is
\begin{equation}
  K(x,y)=\sum_{z=0}^{\min(x,y)}\pi(x-z,N-z,b,c)\pi(z,y,a,b)
  \ \ \Bigl(\Rightarrow\det K=\prod_{n\in\cX}
  \frac{b^n(a\,;q)_n(c\,;q)_n}{(ab\,;q)_n(bc\,;q)_n}\Bigr).
  \label{qh13K}
\end{equation}
For the following $\bm{\lambda}$,
\begin{equation}
  \bm{\lambda}=(ab,c),
  \label{qH1la}
\end{equation}
the symmetry condition \eqref{Kcond1} is satisfied and
{\bf Theorem\,\ref{theo1}} gives \eqref{Keig},
\begin{equation*}
  \sum_{y\in\cX}K(x,y)\pi(y,N,ab,c)\check{P}_n(y,ab,c)
  =\kappa(n)\pi(x,N,ab,c)\check{P}_n(x,ab,c).
\end{equation*}
By evaluating the eigenvalue formula \eqref{eigform1}, we obtain a balanced
${}_3\phi_2$
\begin{equation}
  \kappa(n)=\sum_{y\in\cX}\pi(y,N,b,c)\check{P}_n(y,ab,c)
  ={}_3\phi_2\Bigl(
  \genfrac{}{}{0pt}{}{q^{-n},\,abcq^{n-1},\,b}{ab,\,bc}\Bigm|q\,;q\Bigr)
  =\frac{b^n(a\,;q)_n(c\,;q)_n}{(ab\,;q)_n(bc\,;q)_n}.
  \label{1qHphi}
\end{equation}
The last equality is due to the $q$-Pfaff-Saalsch\"utz formula
\eqref{KLS(1.11.9)}.

\paragraph{$q$-Hahn\,$\to$\,$q$-Meixner}

This limit is achieved by fixing $a$, $b$ with $c\to-c^{-1}q^{1-N}$,
$N\to\infty$, ($c>0$),
\begin{align}
  &\check{P}_n(x,ab,c)\to
  \check{P}_{\text{$q$M}\,n}\bigl(x,abq^{-1},(ab)^{-1}c\bigr),\n
  &\pi(x,N,ab,c)\to\pi_{\text{$q$M}}\bigl(x,abq^{-1},(ab)^{-1}c\bigr),\quad
  \kappa(n)\to\kappa_{\text{$q$M}}(n)=\frac{(a\,;q)_n}{(ab\,;q)_n},
  \label{piqM}\\
  &\pi(x-z,N-z,b,c)\n
  &\qquad\to\pi'_{\text{$q$M}}(x,z,bq^{-1},b^{-1}c)
  =\frac{(-b^{-1}c\,;q)_z(b\,;q)_{x-z}(b^{-1}c)^{x-z}\,
  q^{\binom{x}{2}-\binom{z}{2}}}{(-c\,;q)_x(q\,;q)_{x-z}}
  \frac{(-c\,;q)_\infty}{(-b^{-1}c\,;q)_\infty},
  \label{piqMr1}\\[2pt]
  &K(x,y)\to K_{\text{$q$M}}(x,y,a,b,c)=\sum_{z=0}^{\min(x,y)}
  \pi'_{\text{$q$M}}(x,z,bq^{-1},b^{-1}c)\pi(z,y,a,b).
  \label{KqMr1}
\end{align}
The above $K_{\text{$q$M}}$ \eqref{KqMr1} is not a standard convolution of
$\pi$'s of ($q$M) and ($q$H), as the form of $\pi'_{\text{$q$M}}$ \eqref{piqMr1}
is markedly different from that of \eqref{piqM}.
The relations \eqref{Kcond1}, \eqref{Keig} and \eqref{eigform} of ($q$H)
reduce to those of ($q$M).
Based on $K_{\text{$q$M}}$ \eqref{KqMr1}, let us rederive these results.
The symmetry condition \eqref{Kcond1inf} is satisfied for
$\bm{\lambda}=(abq^{-1},(ab)^{-1}c)$ and {\bf Theorem\,\ref{theo1}} gives
\eqref{Keiginf}.
The eigenvalue formula \eqref{eigforminf} takes a neat form and the sum is
directly evaluated 
\begin{align}
  \kappa_{\text{$q$M}}(n)&=\sum_{y=0}^{\infty}
  \pi_{\text{$q$M}}(y,bq^{-1},b^{-1}c)
  \check{P}_{\text{$q$M}\,n}\bigl(y,abq^{-1},(ab)^{-1}c\bigr)\n
  &={}_2\phi_1\Bigl(
  \genfrac{}{}{0pt}{}{q^{-n},\,b}{ab}\Bigm|q\,;aq^n\Bigr)
  =\frac{(a\,;q)_n}{(ab\,;q)_n}.
\end{align}
In the last equality \eqref{KLS(1.11.4)} is used.
The ${}_2\phi_1$ form is also obtained from the balanced ${}_3\phi_2$ form
\eqref{1qHphi} by the above limit.

The other set of eigenvectors is
\begin{align}
  &\quad(-1)^x\sqrt{\pi_{\text{$q$M}}(x,abq^{-1},(ab)^{-1}c)
  \pi_{\text{$q$M}}(x,-cq^{-1},abc^{-1})}\,
  \check{P}_{\text{$q$M}\,n}(x,-cq^{-1},abc^{-1})\n
  &=(-1)^x\frac{q^{\binom{x}{2}}}{(q;q)_x}
  \check{P}_{\text{$q$M}\,n}(x,-cq^{-1},abc^{-1})\times
  \sqrt{\frac{(ab,-c\,;q)_{\infty}}{(-abc^{-1},-(ab)^{-1}c\,;q)_{\infty}}},
  \label{qMev1}
\end{align}
and the corresponding eigenvalue formula \eqref{eigforminf} reads
\begin{equation}
  \kappa^{(-)}_{\text{$q$M}}(n)=
  \sum_{y=0}^{\infty}K_{\text{$q$M}}(0,y)(-1)^y\sqrt{
  \frac{\pi_{\text{$q$M}}(y,\bm{\lambda})\pi^{(-)}_{\text{$q$M}}(y,\bm{\lambda})}
  {\pi_{\text{$q$M}}(0,\bm{\lambda})\pi^{(-)}_{\text{$q$M}}(0,\bm{\lambda})}}\,
  \check{P}^{(-)}_{\text{$q$M}\,n}(y,\bm{\lambda}).
  \label{eigforminfqM}
\end{equation}
After a few lines of direct calculation, we obtain
\begin{equation}
  \kappa^{(-)}_{\text{$q$M}}(n)
  =\frac{(a,-c\,;q)_{\infty}}{(ab,-b^{-1}c\,;q)_{\infty}}
  {}_2\phi_1\Bigl(\genfrac{}{}{0pt}{}{q^{-n},\,b}{-c}
  \Bigm|q\,;-b^{-1}cq^n\Bigr)
  =\frac{(a,-c\,;q)_{\infty}}{(ab,-b^{-1}c\,;q)_{\infty}}
  \frac{(-b^{-1}c\,;q)_n}{(-c\,;q)_n},
\end{equation}
with
\begin{equation*}
  0<\kappa^{(-)}_{\text{$q$M}}(n)<1,\quad 
  0<\frac{\kappa^{(-)}_{\text{$q$M}}(n)}{\kappa_{\text{$q$M}}(n)}
  =\frac{(aq^n,-cq^n\,;q)_{\infty}}{(abq^n,-b^{-1}cq^n\,;q)_{\infty}}<1.
\end{equation*}

The reduction to $q$-Charlier is not feasible, as it requires $ab\to0$ in
$\check{P}_n(y,ab,c)$, leading $(ab)^{-1}c$ to diverge.

\subsection{Type (\romannumeral2) convolution}

Since this convolution
\begin{equation*}
  \text{(\romannumeral2)}:\ \ K(x,y)\eqdef\sum_{z=\max(0,x+y-N)}^{\min(x,y)}
  \!\!\!\!\!\!\!\!\!\pi(x-z,N-y,\bm{\lambda}_2)\pi(z,y,\bm{\lambda}_1)
\end{equation*}
is not of a factorised form, the determinant of $K$ is not obtained easily.
For $x=0$ only one term $z=0$ contributes and the general eigenvalue formula
\eqref{eigform} takes a simple form
\begin{equation}
  \kappa(n)=\sum_{y\in\cX}\pi(0,N-y,\bm{\lambda}_2)\pi(0,y,\bm{\lambda}_1)
  \frac{\pi(y,N,\bm{\lambda})}{\pi(0,N,\bm{\lambda})}\check{P}_n(y,\bm{\lambda})
  \ \ (n\in\cX).
  \label{eigform2}
\end{equation}
This convolution was employed for Markov chains related with the Hahn polynomial
in \cite{diaconis20}, but the use of the convolution was not stated explicitly.

\subsubsection{Krawtchouk (K)}

By taking $\bm{\lambda}_1=a$ and $\bm{\lambda}_2=b$, the matrix $K(x,y)$ is
\begin{equation}
  K(x,y)=\sum_{z=\max(0,x+y-N)}^{\min(x,y)}\pi(x-z,N-y,b)\pi(z,y,a).
\end{equation}
For the following $\bm{\lambda}$,
\begin{equation}
  \bm{\lambda}=p\eqdef\frac{b}{1-a+b},
\end{equation}
the symmetry condition \eqref{Kcond1} is satisfied and
{\bf Theorem\,\ref{theo1}} gives \eqref{Keig}.
By directly evaluating the eigenvalue formula \eqref{eigform2}, we obtain
\begin{equation}
  \kappa(n)=\sum_{y\in\cX}\pi(y,N,b)\check{P}_n(y,p)=(a-b)^n
  ={}_1F_0\Bigl(\genfrac{}{}{0pt}{}{-n}{-}\Bigm|bp^{-1}\Bigr).
\end{equation}
The use of the generating function \eqref{KLS(9.11.11)} and the self-duality (K)
leads to the above result, too.
It is interesting to note that odd eigenvalues are all negative if $0<a<b<1$.

\paragraph{Krawtchouk\,$\to$\,Charlier}

This is exactly the same situation of (K)$\to$(C) \eqref{KCconv1} in the type
(\romannumeral1) convolution.

\subsubsection{Hahn (H)}

By taking $\bm{\lambda}_1=(a,b)$ and $\bm{\lambda}_2=(b,c)$, the matrix
$K(x,y)$ is
\begin{equation}
  K(x,y)=\!\!\!\sum_{z=\max(0,x+y-N)}^{\min(x,y)}\!\!\!
  \pi(x-z,N-y,b,c)\pi(z,y,a,b).
\end{equation}
The special cases of this convolution with $a=b=c=\frac12$ corresponding to
the discrete Chebyshev polynomials and $a=b=c=\frac12\theta$ are discussed in
detail in \cite{diaconis20}.
For the following $\bm{\lambda}$,
\begin{equation}
  \bm{\lambda}=(a+b,b+c),
\end{equation}
the symmetry condition \eqref{Kcond1} is satisfied and
{\bf Theorem\,\ref{theo1}} gives \eqref{Keig}.
By writing down the eigenvalue formula \eqref{eigform2}, we have
\begin{align}
  \kappa(n)&=\sum_{y\in\cX}\pi(y,N,b,c)\check{P}_n(y,a+b,b+c)
  ={}_3F_2\Bigl(
  \genfrac{}{}{0pt}{}{-n,\,n+a+2b+c-1,\,b}{a+b,\,b+c}\Bigm|1\Bigr)\\
  &=\sum_{k=0}^n\binom{n}{k}(-1)^k\frac{(b)_k(n+a+2b+c-1)_k}{(a+b)_k(b+c)_k}
  \label{gendiaco}.
\end{align}
The above explicit eigenvalue formula \eqref{gendiaco} reproduces the known
results (5.9) in \cite{diaconis20},
which corresponds to the special case $a=b=c=\frac12\theta$.

\paragraph{Hahn\,$\to$\,Meixner}

By fixing $a$, $b$ with $c\to N(1-c)c^{-1}$ ($\Rightarrow 0<c<1$) and taking
the limit $N\to\infty$,
one obtains the same Meixner limit $K_{\text{M}}(x,y)$ as in \eqref{HM1}
\begin{equation*}
  K(x,y)\to K_{\text{M}}(x,y,a,b,c)
  =\sum_{z=0}^{\min(x,y)}\pi_{\text{M}}(x-z,b,c)\pi(z,y,a,b).
\end{equation*}

It is interesting to note that the $q$-Hahn version of this convolution does
not work, due to the presence of the extra factor $a^{N-x}$ in $\pi(x,N,a,b)$
in \eqref{qHpi}.

\subsection{Type (\romannumeral3) convolution}

Since this convolution
\begin{equation*}
  \text{(\romannumeral3)}:\ \ K(x,y)\eqdef\sum_{z=\max(x,y)}^N
  \!\!\!\!\pi(x,z,\bm{\lambda}_2)\pi(z-y,N-y,\bm{\lambda}_1),
\end{equation*}
is of a factorised form, the first factor $\pi(x,z,\bm{\lambda}_2)$ being upper
triangular and the second factor $\pi(z-y,N-y,\bm{\lambda}_1)$ lower triangular,
the eigenvalues can be guessed from the determinant,
\begin{equation*}
  \prod_{n\in\cX}\kappa(n)=\det K
  =\prod_{n\in\cX}\pi(n,n,\bm{\lambda}_2)\pi(0,n,\bm{\lambda}_1).
\end{equation*}

\subsubsection{Krawtchouk (K)}

By taking $\bm{\lambda}_1=a$ and $\bm{\lambda}_2=b$, the matrix $K(x,y)$ is
\begin{align}
  K(x,y)&=\sum_{z=\max(x,y)}^N\pi(x,z,b)\pi(z-y,N-y,a)\quad
  \Bigl(\Rightarrow\det K=\prod_{n\in\cX}(1-a)^nb^n\Bigr).
  \label{KK3}
\end{align}
For the following $\bm{\lambda}$,
\begin{equation}
  \bm{\lambda}=p\eqdef\frac{ab}{1-b+ab},
\end{equation}
the symmetry condition \eqref{Kcond1} is satisfied and
{\bf Theorem\,\ref{theo1}} gives \eqref{Keig}.
By evaluating the eigenvalue formula \eqref{eigform}, we have
\begin{equation}
  \kappa(n)=\sum_{z=0}^N\pi(z,N,a)\sum_{y=0}^z\pi(y,z,b)\check{P}_n(y,p)
  =(1-a)^nb^n
  = {}_1F_0\Bigl(\genfrac{}{}{0pt}{}{-n}{-}\Bigm|abp^{-1}\Bigr).
\end{equation}
The use of the generating function
\eqref{KLS(9.11.11)} and the self-duality (K) give the same result.

\paragraph{Krawtchouk\,$\to$\,Charlier}

This is achieved by $a\to aN^{-1}$, $N\to\infty$,
\begin{align}
  &\check{P}_n(x,p)\to P_{\text{C}\,n}(x,p'),\quad
  p'\eqdef\frac{ab}{1-b},\n
  &\pi(x,N,p)\to\pi_{\text{C}}(x,p'),\quad
  \kappa(n)\to\kappa_{\text{C}}(n)=b^n,\n
  &K(x,y)\to K_{\text{C}}(x,y,a,b)
  =\sum_{z=\max(x,y)}^{\infty}\pi(x,z,b)\pi_{\text{C}}(z-y,a),
  \label{KCconv3}
\end{align}
and the relations \eqref{Kcond1}, \eqref{Keig} and \eqref{eigform} of (K)
reduce to those of (C).
Based on $K_{\text{C}}$ \eqref{KCconv3}, let us rederive these results.
The symmetry condition \eqref{Kcond1inf} is satisfied for $\bm{\lambda}=p'$ and
{\bf Theorem\,\ref{theo1}} gives \eqref{Keiginf}.
The eigenvalue formula \eqref{eigforminf} reads
\begin{equation}
  \kappa_{\text{C}}(n)
  =\sum_{z=0}^{\infty}\pi_{\text{C}}(z,a)\sum_{y=0}^z\pi(y,z,b)
  \check{P}_{\text{C}\,n}(y,p')=b^n
  ={}_1F_0\Bigl(\genfrac{}{}{0pt}{}{-n}{-}\Bigm|abp^{\prime\,-1}\Bigr),
\end{equation}
by using the generating function \cite{kls}(9.14.11) and the self-duality (C).

\subsubsection{Hahn (H)}

By taking $\bm{\lambda}_1=(a,b)$ and $\bm{\lambda}_2=(c,a)$, the matrix
$K(x,y)$ is
\begin{equation}
  K(x,y)=\!\!\!\sum_{z=\max(x,y)}^N\!\!\!\pi(x,z,c,a)\pi(z-y,N-y,a,b)
  \ \ \Bigl(\Rightarrow\det K=\prod_{n\in\cX}\frac{(b)_n(c)_n}{(a+b)_n(a+c)_n}
  \Bigr).
\end{equation}
For the following $\bm{\lambda}$,
\begin{equation}
  \bm{\lambda}=(c,a+b),
\end{equation}
the symmetry condition \eqref{Kcond1} is satisfied and
{\bf Theorem\,\ref{theo1}} gives \eqref{Keig}.
By evaluating the eigenvalue formula \eqref{eigform}, we obtain a balanced
${}_3F_2$
\begin{align}
  \kappa(n)&=\sum_{z=0}^N\pi(z,N,a,b)\sum_{y=0}^z\pi(y,z,c,a)
  \check{P}_n(y,c,a+b)\n
  &={}_3F_2\Bigl(\genfrac{}{}{0pt}{}{-n,\,n+a+b+c-1,\,a}{a+b,\,a+c}\Bigm|1\Bigr)
  =\frac{(b)_n(c)_n}{(a+b)_n(a+c)_n}.
  \label{3HF}
\end{align}
In the last equality Pfaff-Saalsch\"utz formula \eqref{KLS(1.5.5)} is used.

\paragraph{Hahn\,$\to$\,Meixner}

This is achieved by fixing $a$ and $c$ with $b\to N(1-b)b^{-1}$
($\Rightarrow 0<b<1$), $N\to\infty$,
\begin{align}
  &\check{P}_n(x,c,a+b)\to\check{P}_{\text{M}\,n}(x,c,b),\n
  &\pi(x,N,c,a+b)\to\pi_{\text{M}}(x,c,b),\quad
  \kappa(n)\to\kappa_{\text{M}}(n)=\frac{(c)_n}{(a+c)_n},\n
  &K(x,y)\to K_{\text{M}}(x,y,a,b,c)=\sum_{z=\max(x,y)}^{\infty}\pi(x,z,c,a)
  \pi_{\text{M}}(z-y,a,b),
  \label{HM3}
\end{align}
and the relations \eqref{Kcond1}, \eqref{Keig} and \eqref{eigform} of (H)
reduce to those of (M).
This is a standard convolution of $\pi$'s of (H) and (M), but the order is
opposite from that of \eqref{HM1}.
Based on $K_{\text{M}}$ \eqref{HM3}, let us rederive these results.
The symmetry condition \eqref{Kcond1inf} is satisfied for
$\bm{\lambda}=(c,b)$ and {\bf Theorem\,\ref{theo1}} gives \eqref{Keiginf}.
The eigenvalue formula \eqref{eigforminf} reads
\begin{equation}
  \kappa_{\text{M}}(n)=\sum_{z=0}^{\infty}\pi_{\text{M}}(z,a,b)\sum_{y=0}^z
  \pi(y,z,c,a)\check{P}_{\text{M}\,n}(y,c,b)
  ={}_2F_1\Bigl(\genfrac{}{}{0pt}{}{-n,\,a}{a+c}\Bigm|1\Bigr)
  =\frac{(c)_n}{(a+c)_n}.
\end{equation}
Of course the ${}_2F_1$ form is also obtained from the balanced ${}_3F_2$ form
\eqref{3HF} by the above limit.

\paragraph{Hahn\,$\to$\,Meixner\,$\to$\,Charlier}

This is achieved by $a\to aN$, $c\to cN$, $b\to\frac{b}{b+N}$, $N\to\infty$,
\begin{align}
  &P_{\text{M}\,n}\bigl(x,cN,\tfrac{b}{b+N}\bigr)\to P_{\text{C}\,n}(x,bc),\n
  &\pi_{\text{M}}\bigl(x,cN,\tfrac{b}{b+N}\bigr)\to\pi_{\text{C}}(x,bc),\quad
  \kappa_{\text{M}}(n)=\frac{(cN)_n}{\bigl((a+c)N\bigr)_n}\to
  \kappa_{\text{C}}(n)=\Bigl(\frac{c}{a+c}\Bigr)^n,\n
  &\pi(x,z,cN,aN)\to\pi_{\text{K}}\bigl(x,z,\tfrac{c}{a+c}\bigr),\quad
  \pi_{\text{M}}(x,aN,\tfrac{b}{b+N}\bigr)
  \to\pi_{\text{C}}(x,ab),\n
  &K_{\text{M}}(x,y,a,b,c)\to K_{\text{C}}(x,y,a,b,c)
  =\sum_{z=\max(x,y)}^{\infty}\pi_{\text{K}}(x,z,p)
  \pi_{\text{C}}(z-y,ab),
  \label{HMC2}
\end{align}
and the relations \eqref{Kcond1inf}, \eqref{Keiginf} and \eqref{eigforminf}
of (M) reduce to those of (C).
This $K_{\text{C}}$ agrees with \eqref{KCconv3} with the replacement
$(a,b)\to(ab,p)$.

\subsubsection{$q$-Hahn ($q$H)}
\label{3qH}

By taking $\bm{\lambda}_1=(a,b)$ and $\bm{\lambda}_2=(c,a)$, the matrix
$K(x,y)$ is
\begin{equation} 
  K(x,y)=\!\!\!\sum_{z=\max(x,y)}^N\!\!\!\pi(x,z,c,a)\pi(z-y,N-y,a,b)
  \ \ \Bigl(\Rightarrow\det K=\prod_{n\in\cX}
  \frac{a^n(b\,;q)_n(c\,;q)_n}{(ab\,;q)_n(ac\,;q)_n}\Bigr).
\end{equation}
For the following $\bm{\lambda}$,
\begin{equation}
  \bm{\lambda}=(c,ab),
\end{equation}
the symmetry condition \eqref{Kcond1} is satisfied and
{\bf Theorem\,\ref{theo1}} gives \eqref{Keig}.
By evaluating the eigenvalue formula \eqref{eigform}, we obtain a balanced
${}_3\phi_2$
\begin{align}
  \kappa(n)&=\sum_{z=0}^N\pi(z,N,a,b)\sum_{y=0}^z\pi(y,z,c,a)
  \check{P}_n(y,c,ab)\n
  &={}_3\phi_2\Bigl(\genfrac{}{}{0pt}{}{q^{-n},\,abcq^{n-1},\,a}{ac,\,ab}
  \Bigm|q\,;q\Bigr)
  =\frac{a^n(b\,;q)_n(c\,;q)_n}{(ab\,;q)_n(ac\,;q)_n}.
  \label{3qHphi}
\end{align}
The last equality is due to the $q$-Pfaff-Saalsch\"utz formula
\eqref{KLS(1.11.9)}.

\paragraph{$q$-Hahn\,$\to$\,$q$-Meixner}

This limit is achieved by fixing $a$, $c$ with $b\to-b^{-1}q^{1-N}$,
$N\to\infty$, ($b>0$),
\begin{align}
  &\check{P}_n(x,c,ab)\to
  \check{P}_{\text{$q$M}\,n}\bigl(x,cq^{-1},(ac)^{-1}b\bigr),\n
  &\pi(x,N,c,ab)\to\pi_{\text{$q$M}}\bigl(x,cq^{-1},(ac)^{-1}b\bigr),\quad
  \kappa(n)\to\kappa_{\text{$q$M}}(n)=\frac{(c\,;q)_n}{(ac\,;q)_n},\n
  &\pi(z-y,N-y,a,b)\n
  &\qquad\to\pi'_{\text{$q$M}}(z,y,aq^{-1},a^{-1}b)
  =\frac{(-a^{-1}b\,;q)_y(a\,;q)_{z-y}(a^{-1}b)^{z-y}
  q^{\binom{z}{2}-\binom{y}{2}}}{(-b\,;q)_z(q\,;q)_{z-y}}
  \frac{(-b\,;q)_{\infty}}{(-a^{-1}b\,;q)_{\infty}},\n
  &K(x,y)\to K_{\text{$q$M}}(x,y,a,b,c)=\!\!\!\sum_{z=\max(x,y)}^{\infty}\!\!\!
  \pi(x,z,c,a)\pi'_{\text{$q$M}}(z,y,aq^{-1},a^{-1}b).
  \label{KqMr2}
\end{align}
This is not a standard convolution as $\pi'_{\text{$q$M}}$ is not an
orthogonality measure of ($q$M).
The relations \eqref{Kcond1}, \eqref{Keig} and \eqref{eigform} of ($q$H)
reduce to those of ($q$M).
Based on $K_{\text{$q$M}}$ \eqref{KqMr2}, let us rederive these results.
The symmetry condition \eqref{Kcond1inf} is satisfied for
$\bm{\lambda}=(cq^{-1},(ac)^{-1}b)$ and
{\bf Theorem\,\ref{theo1}} gives \eqref{Keiginf}.
The eigenvalue formula \eqref{eigforminf} is written down as
\begin{align}
  \kappa_{\text{$q$M}}(n)&=\sum_{z=0}^{\infty}
  \pi_{\text{$q$M}}(z,aq^{-1},a^{-1}b)\sum_{y=0}^z\pi(y,z,c,a)
  \check{P}_{\text{$q$M}\,n}\bigl(y,cq^{-1},(ac)^{-1}b\bigr)
  \label{(iii)kappaqM}\n
  &={}_2\phi_1\Bigl(\genfrac{}{}{0pt}{}{q^{-n},\,a}{ac}\Bigm|q\,;cq^n\Bigr)
  =\frac{(c\,;q)_n}{(ac\,;q)_n}.
\end{align}
In the last equality \eqref{KLS(1.11.4)} is used.
The ${}_2\phi_1$ form is also obtained from the balanced ${}_3\phi_2$ form
\eqref{3qHphi} by the above limit.

The other set of eigenvectors of $K_{\text{$q$M}}$ is
\begin{align}
  &\quad(-1)^x\sqrt{\pi_{\text{$q$M}}(x,cq^{-1},(ac)^{-1}b)
  \pi_{\text{$q$M}}(x,-a^{-1}bq^{-1},acb^{-1})}\,
  \check{P}_{\text{$q$M}\,n}(x,-a^{-1}bq^{-1},acb^{-1})\n
  &=(-1)^x\frac{q^{\binom{x}{2}}}{(q;q)_x}
  \check{P}_{\text{$q$M}\,n}(x,-a^{-1}bq^{-1},acb^{-1})\times
  \sqrt{\frac{(-a^{-1}b,c\,;q)_{\infty}}{(-acb^{-1},-(ac)^{-1}b\,;q)_{\infty}}},
\end{align}
and the corresponding eigenvalue formula \eqref{eigforminf} takes the same form
as \eqref{eigforminfqM}.
After a few lines of direct calculation, we obtain
\begin{equation}
  \kappa^{(-)}_{\text{$q$M}}(n)
  =\frac{(c,-b\,;q)_{\infty}}{(ac,-a^{-1}b\,;q)_{\infty}}
  {}_2\phi_1\Bigl(\genfrac{}{}{0pt}{}{q^{-n},\,a}{-b}
  \Bigm|q\,;-a^{-1}bq^n\Bigr)
  =\frac{(c,-b\,;q)_{\infty}}{(ac,-a^{-1}b\,;q)_{\infty}}
  \frac{(-a^{-1}b\,;q)_n}{(-b\,;q)_n},\!
\end{equation}
with
\begin{equation*}
  0<\kappa^{(-)}_{\text{$q$M}}(n)<1,\quad 
  0<\frac{\kappa^{(-)}_{\text{$q$M}}(n)}{\kappa_{\text{$q$M}}(n)}
  =\frac{(cq^n,-bq^n\,;q)_{\infty}}{(acq^n,-a^{-1}bq^n\,;q)_{\infty}}<1.
\end{equation*}
The $q$-Charlier limit does not exist, as it requires $c\to0$.
This causes $\kappa_{\text{$q$M}}(n)\to1$ and $\pi(x,z,c,a)\to0$ in
$K_{\text{$q$M}}$.

\subsection{Type (\romannumeral4) convolution}

This type of convolutions
\begin{equation*}
  \text{(\romannumeral4)}:\ \ K(x,y)\eqdef\sum_{z_2=0}^{\min(x,y)}
  \!\!\pi(z_2,y,\bm{\lambda}_1)\!\!\!\sum_{z_1=\max(x,y)}^N
  \!\!\!\!\pi(x-z_2,z_1-z_2,\bm{\lambda}_3)\pi(z_1-y,N-y,\bm{\lambda}_2),
\end{equation*}
was reported in \cite{hoa-rah83} for (K).
The eigenvalues were derived by a different method from that given below.
At $x=0$ only the $z_2=0$ term contributes and the eigenvalue formula becomes a
manageable double sum formula
\begin{equation}
  \kappa(n)=\sum_{z=0}^N\sum_{y=0}^z\pi(0,y,\bm{\lambda}_1)
  \pi(0,z,\bm{\lambda}_3)\pi(z-y,N-y,\bm{\lambda}_2)
  \frac{\pi(y,N,\bm{\lambda})}{\pi(0,N,\bm{\lambda})}
  \check{P}_n(y,\bm{\lambda}).
  \label{eigform4}
\end{equation}

\subsubsection{Krawtchouk (K)}

By taking $\bm{\lambda}_1=a$, $\bm{\lambda}_2=b$ and $\bm{\lambda}_3=c$,
the matrix $K(x,y)$ is
\begin{align}
  K(x,y)&=\sum_{z_2=0}^{\min(x,y)}\pi(z_2,y,a)
  \sum_{z_1=\max(x,y)}^N\pi(x-z_2,z_1-z_2,c)\pi(z_1-y,N-y,b).
\end{align}
For the following $\bm{\lambda}$,
\begin{equation}
  \bm{\lambda}=p\eqdef\frac{bc}{bc+(1-a)(1-c)},
\end{equation}
the symmetry condition \eqref{Kcond1} is satisfied and
{\bf Theorem\,\ref{theo1}} gives \eqref{Keig}.
By writing down the eigenvalue formula \eqref{eigform4}, we have
\begin{align}
  \kappa(n)&=\sum_{z=0}^N\pi(z,N,b)\sum_{y=0}^z\pi(y,z,c)\check{P}_n(y,p)
  ={}_1F_0\Bigl(\genfrac{}{}{0pt}{}{-n}{-}\Bigm|bcp^{-1}\Bigr)\n
  &=(a+c-ac-bc)^n=\bigl(1-bc-(1-a)(1-c)\bigr)^n,
\end{align}
by the generating function \eqref{KLS(9.11.11)} and the self-duality (K).

\paragraph{Krawtchouk\,$\to$\,Charlier}

This is achieved by fixing $a$ and $c$ with $b\to bN^{-1}$, $N\to\infty$,
\begin{align}
  &\check{P}_n(x,p)\to\check{P}_{\text{C}\,n}(x,p'),\quad
  p'\eqdef\frac{bc}{(1-a)(1-c)},\n
  &\pi(x,N,p)\to\pi_{\text{C}}(x,p'),\quad
  \kappa(n)\to\kappa_{\text{C}}(n)=(a+c-ac)^n=\bigl(1-(1-a)(1-c)\bigr)^n,\n
  &K(x,y)\to K_{\text{C}}(x,y,a,b,c)=\!\!\!\sum_{z_2=0}^{\min(x,y)}\!\!\!
  \pi(z_2,y,a)\!\!\!\sum_{z_1=\max(x,y)}^{\infty}\!\!\!\!
  \pi(x-z_2,z_1-z_2,c)\pi_{\text{C}}(z_1-y,b),
  \label{KCconv4}
\end{align}
and the relations \eqref{Kcond1}, \eqref{Keig} and \eqref{eigform} of (K)
reduce to those of (C).
Based on $K_{\text{C}}$ \eqref{KCconv4}, let us rederive these results.
The symmetry condition \eqref{Kcond1inf} is satisfied for $\bm{\lambda}=p'$ and
{\bf Theorem\,\ref{theo1}} gives \eqref{Keiginf}.
The eigenvalue formula \eqref{eigforminf} reads
\begin{equation}
  \kappa_{\text{C}}(n)
  =\sum_{z=0}^{\infty}\pi_{\text{C}}(z,b)\sum_{y=0}^z\pi(y,z,c)
  \check{P}_{\text{C}\,n}(y,p')=(a+c-ac)^n=
  {}_1F_0\Bigl(\genfrac{}{}{0pt}{}{-n}{-}\Bigm|bcp^{\prime\,-1}\Bigr),
\end{equation}
by using the generating function \cite{kls}(9.14.11) and the self-duality (C).

\subsubsection{Hahn (H)}

By taking $\bm{\lambda}_1=(a_1,b_1)$, $\bm{\lambda}_2=(a_2,b_2)$ and
$\bm{\lambda}_3=(b_1,a_2)$, the matrix $K(x,y)$ with four parameters is
\begin{equation}
  K(x,y)=\!\!\!\sum_{z_2=0}^{\min(x,y)}\!\!\!\pi(z_2,y,a_1,b_1)
  \!\!\!\sum_{z_1=\max(x,y)}^N\!\!\!\pi(x-z_2,z_1-z_2,b_1,a_2)
  \pi(z_1-y,N-y,a_2,b_2).
\end{equation}
For the following $\bm{\lambda}$,
\begin{equation}
  \bm{\lambda}=(a_1+b_1,a_2+b_2),
\end{equation}
the symmetry condition \eqref{Kcond1} is satisfied and
{\bf Theorem\,\ref{theo1}} gives \eqref{Keig}.
By evaluating the eigenvalue formula \eqref{eigform4}, we obtain a balanced
${}_4F_3$
\begin{align}
  \kappa(n)&=\sum_{z=0}^N\pi(z,N,a_2,b_2)\sum_{y=0}^z\pi(y,z,b_1,a_2)
  \check{P}_n(y,a_1+b_1,a_2+b_2)\n
  &={}_4F_3\Bigl(\genfrac{}{}{0pt}{}{-n,\,n+a_1+b_1+a_2+b_2-1,\,b_1,\,a_2}
  {a_1+b_1,\,b_1+a_2,\,a_2+b_2}\Bigm|1\Bigr).
\end{align}

\paragraph{Hahn\,$\to$\,Meixner}

This is achieved by fixing $a_1$, $b_1$ and $a_2$ with $b_2\to N(1-b_2)b_2^{-1}$
($\Rightarrow 0<b_2<1$), $N\to\infty$,
\begin{align}
  &\check{P}_n(x,a_1+b_1,a_2+b_2)\to\check{P}_{\text{M}\,n}(x,a_1+b_1,b_2),\n
  &\pi(x,N,a_1+b_1,a_2+b_2)\to\pi_{\text{M}}(x,a_1+b_1,b_2),\n
  &K(x,y)\to K_{\text{M}}(x,y,a_1,b_1,a_2,b_2)\n
  &\phantom{K(x,y)\to}
  =\sum_{z_2=0}^{\min(x,y)}\!\!\pi(z_2,y,a_1,b_1)
  \!\!\!\sum_{z_1=\max(x,y)}^{\infty}\!\!\!\!\pi(x-z_2,z_1-z_2,b_1,a_2)
  \pi_{\text{M}}(z_1-y,a_2,b_2),
  \label{HM4}
\end{align}
and the relations \eqref{Kcond1}, \eqref{Keig} and \eqref{eigform} of (H)
reduce to those of (M).
Based on $K_{\text{M}}$ \eqref{HM4}, let us rederive these results.
The symmetry condition \eqref{Kcond1inf} is satisfied for
$\bm{\lambda}=(a_1+b_1,b_2)$ and {\bf Theorem\,\ref{theo1}} gives
\eqref{Keiginf}.
The eigenvalue formula \eqref{eigforminf} is written down as
\begin{align}
  \kappa_{\text{M}}(n)&=\sum_{z=0}^{\infty}\pi_{\text{M}}(z,a_2,b_2)\sum_{y=0}^z
  \pi(y,z,b_1,a_2)\check{P}_{\text{M}\,n}(y,a_1+b_1,b_2)
  \label{(iv)kappaM}\n
  &={}_3F_2\Bigl(\genfrac{}{}{0pt}{}{-n,\,b_1,\,a_2}
  {a_1+b_1,\,b_1+a_2}\Bigm|1\Bigr).
\end{align}

\paragraph{Hahn\,$\to$\,Meixner\,$\to$\,Charlier}

This is achieved by $a_1\to a_1N$, $b_1\to b_1N$, $a_2\to a_2N$,
$b_2\to\frac{b_2}{b_2+N}$, $N\to\infty$,
\begin{align}
  &\check{P}_{\text{M}\,n}\bigl(x,(a_1+b_1)N,\tfrac{b_2}{b_2+N}\bigr)\to
  P_{\text{C}\,n}\bigl(x,(a_1+b_1)b_2\bigr),\n
  &\pi_{\text{M}}\bigl(x,a_2N,\tfrac{b_2}{b_2+N}\bigr)
  \to\pi_{\text{C}}(x,a_2b_2),\n
  &\pi(z,y,a_1N,b_1N)\to\pi_{\text{K}}\bigl(z,y,\tfrac{a_1}{a_1+b_1}\bigr),\quad
  \pi(x,z,b_1N,a_2N)\to\pi_{\text{K}}\bigl(x,z,\tfrac{b_1}{b1+a_2}\bigr),\n
  &K_{\text{M}}(x,y)\to K_{\text{C}}(x,y)
  \label{HMC4}\\
  &\phantom{K_{\text{M}}(x,y)\to}
  =\sum_{z_2=0}^{\min(x,y)}\!\!
  \pi_{\text{K}}\bigl(z_2,y,\tfrac{a_1}{a_1+b_1}\bigr)
  \!\!\sum_{z_1=\max(x,y)}^{\infty}\!\!\!\!\!
  \pi_{\text{K}}\bigl(x-z_2,z_1-z_2,\tfrac{b_1}{b_1+a_2}\bigr)
  \pi_{\text{C}}(z_1-y,a_2b_2),
  \nonumber
\end{align}
and the relations \eqref{Kcond1inf}, \eqref{Keiginf} and \eqref{eigforminf}
of (M) reduce to those of (C).
This $K_{\text{C}}$ agrees with \eqref{KCconv4} with the replacement
$(a,b,c)\to(\frac{a_1}{a_1+b_1},a_2b_2,\frac{b_1}{b_1+a_2})$.

\subsubsection{$q$-Hahn ($q$H)}
\label{4qH}

This convolution for $q$-Hahn has almost the same structure as that for Hahn.
By taking $\bm{\lambda}_1=(a_1,b_1)$, $\bm{\lambda}_2=(a_2,b_2)$ and
$\bm{\lambda}_3=(b_1,a_2)$, the matrix $K(x,y)$ with four parameters is
\begin{equation}
  K(x,y)=\!\!\!\sum_{z_2=0}^{\min(x,y)}\!\!\!\pi(z_2,y,a_1,b_1)
  \!\!\!\sum_{z_1=\max(x,y)}^N\!\!\!\pi(x-z_2,z_1-z_2,b_1,a_2)
  \pi(z_1-y,N-y,a_2,b_2).
\end{equation}
For the following $\bm{\lambda}$,
\begin{equation}
  \bm{\lambda}=(a_1b_1,a_2b_2),
\end{equation}
the symmetry condition \eqref{Kcond1} is satisfied and
{\bf Theorem\,\ref{theo1}} gives \eqref{Keig}.
The eigenvalue formula \eqref{eigform4} takes a neat form and after a few
lines of direct calculation, we obtain a balanced ${}_4\phi_3$,
\begin{align}
  \kappa(n)&=\sum_{z=0}^N\pi(z,N,a_2,b_2)\sum_{y=0}^z\pi(y,z,b_1,a_2)
  \check{P}_n(y,a_1b_1,a_2b_2)\n
  &={}_4\phi_3\Bigl(\genfrac{}{}{0pt}{}
  {q^{-n},\,a_1b_1a_2b_2q^{n-1},\,b_1,\,a_2}
  {a_1b_1,\,b_1a_2,\,a_2b_2}\Bigm|q\,;q\Bigr)
  \label{(4.109)}.
\end{align}

\paragraph{$q$-Hahn\,$\to$\,$q$-Meixner}

This limit is achieved by fixing $a_1$, $b_1$, $a_2$ with
$b_2\to-b_2^{-1}q^{1-N}$, $N\to\infty$,
\begin{align}
  &\check{P}_n(x,a_1b_1,a_2b_2)\to
  \check{P}_{\text{$q$M}\,n}\bigl(x,a_1b_1q^{-1},(a_1b_1a_2)^{-1}b_2\bigr),\n
  &\pi(x,N,a_1b_1,a_2b_2)\to
  \pi_{\text{$q$M}}\bigl(x,a_1b_1q^{-1},(a_1b_1a_2)^{-1}b_2\bigr),\n
  &\pi(z_1-y,N-y,a_2,b_2)\n
  &\qquad\to\pi'_{\text{$q$M}}(z_1,y,a_2q^{-1},a_2^{-1}b_2)
  =\frac{(-a_2^{-1}b_2\,;q)_y(a_2\,;q)_{z_1-y}(a_2^{-1}b_2)^{z_1-y}
  q^{\binom{z_1}{2}-\binom{y}{2}}}{(-b_2\,;q)_{z_1}(q\,;q)_{z_1-y}}
  \frac{(-b_2;q)_{\infty}}{(-a_2^{-1}b_2\,;q)_{\infty}},\n
  &K(x,y)\to K_{\text{$q$M}}(x,y)=\!\!\sum_{z_2=0}^{\min(x,y)}\!\!\!
  \pi(z_2,y,a_1,b_1)\!\!\!\!\!\sum_{z_1=\max(x,y)}^{\infty}\!\!\!\!\!
  \pi(x-z_2,z_1-z_2,b_1,a_2)\n
  &\hspace{85mm}\times
  \pi'_{\text{$q$M}}(z_1,y,a_2q^{-1},a_2^{-1}b_2),
  \label{qHqM}
\end{align}
and the relations \eqref{Kcond1}, \eqref{Keig} and \eqref{eigform} of ($q$H)
reduce to those of ($q$M).
The above $K_{\text{$q$M}}$ \eqref{qHqM} is not a convolution of the
orthogonality measures as $\pi'_{\text{$q$M}}$ is not.
Based on $K_{\text{$q$M}}$ \eqref{qHqM}, let us rederive these results.
The symmetry condition \eqref{Kcond1inf} is satisfied for
$\bm{\lambda}=(a_1b_1q^{-1},(a_1b_1a_2)^{-1}b_2)$ and
{\bf Theorem\,\ref{theo1}} gives \eqref{Keiginf}.
The eigenvalue formula \eqref{eigforminf}
takes a neat form and several lines of direct calculation gives
\begin{align}
  \kappa_{\text{$q$M}}(n)&=\sum_{z=0}^{\infty}
  \pi_{\text{$q$M}}(z,a_2q^{-1},a_2^{-1}b_2)\sum_{y=0}^z\pi(y,z,b_1,a_2)
  \check{P}_{\text{$q$M}\,n}\bigl(y,a_1b_1q^{-1},(a_1b_1a_2)^{-1}b_2\bigr)\n
  &={}_3\phi_2\Bigl(\genfrac{}{}{0pt}{}
  {q^{-n},\,b_1,\,a_2}{a_1b_1,\,b_1a_2}\Bigm|q\,;a_1b_1q^n\Bigr).
  \label{(4.114)}
\end{align}
The ${}_3\phi_2$ form \eqref{(4.114)} is also obtained from the ${}_4\phi_3$
form \eqref{(4.109)} by the above mentioned limit.

The other set of eigenvectors of $K_{\text{$q$M}}$ is
\begin{align}
  \!\!\!\!\!&
  \quad(-1)^x\sqrt{\pi_{\text{$q$M}}(x,a_1b_1q^{-1},(a_1b_1a_2)^{-1}b_2)
  \pi_{\text{$q$M}}(x,-a_2^{-1}b_2q^{-1},a_1b_1a_2b_2^{-1})}\n
  \!\!\!\!\!&\hspace{82mm}\times
  \check{P}_{\text{$q$M}\,n}(x,-a_2^{-1}b_2q^{-1},a_1b_1a_2b_2^{-1})\n
  \!\!\!\!\!&=(-1)^x\frac{q^{\binom{x}{2}}}{(q;q)_x}
  \check{P}_{\text{$q$M}\,n}(x,-a_2^{-1}b_2q^{-1},a_1b_1a_2b_2^{-1})
  \times\!
  \sqrt{\frac{(a_1b_1,-a_2^{-1}b_2\,;q)_{\infty}}
  {(-a_1b_1a_2b_2^{-1},-(a_1b_1a_2)^{-1}b_2\,;q)_{\infty}}}.\!
\end{align}
After a few pages of long calculation, we obtain the corresponding eigenvalue
\begin{align}
  \kappa^{(-)}_{\text{$q$M}}(n)&=
  \frac{(b_1,-b_2\,;q)_{\infty}}{(b_1a_2,-a_2^{-1}b_2\,;q)_{\infty}}
  \sum_{z=0}^{\infty}\pi_{\text{$q$M}}(z,a_2q^{-1},-b_1)
  \sum_{y=0}^z\pi(y,z,-a_2^{-1}b_2,a_2)\pi(0,y,a_1,b_1)\n
  &\hspace{80mm}\times
  \check{P}_{\text{$q$M}\,n}(y,-a_2^{-1}b_2q^{-1},a_1b_1a_2b_2^{-1})\n
  &=\frac{(b_1,-b_2\,;q)_{\infty}}{(b_1a_2,-a_2^{-1}b_2\,;q)_{\infty}}
  \sum_{k=0}^n\frac{(q^{-n},b_1,a_2\,;q)_k}{(-b_2,a_1b_1\,;q)_k}
  \frac{(-a_2^{-1}b_2q^n)^k}{(q\,;q)_k}\n
  &\hspace{70mm}\times
  {}_3\phi_2\Bigl(\genfrac{}{}{0pt}{}{a_1,\,a_2q^k,\,-a_2^{-1}b_2q^k}
  {-b_2q^k,\,a_1b_1q^k}\Bigm|q\,;b_1\Bigr).
  \label{viqMka-}
\end{align}
The $q$-Charlier limit does not exist, as it requires $a_1b_1\to0$ in
$\pi_{\text{$q$M}}(x,a_1b_1,a_2^{-1}b_2)$. 
This causes $\pi(z_2,y,a_1,b_1)\pi(x-z_2,z_1-z_2,b_1,a_2)\to0$ in
$K_{\text{$q$M}}$.

\subsection{Type (\romannumeral5) convolution}
\label{convtype5}

We find only one $K(x,y)$ by this convolution for Krawtchouk (K).
By taking $\bm{\lambda}_1=a$, $\bm{\lambda}_2=b$ and $\bm{\lambda}_3=c$,
the matrix $K(x,y)$ is
\begin{equation}
  \text{(\romannumeral5)}:\ \ K(x,y)=\!\!\sum_{z_2=0}^{\min(x,y)}\!\!
  \pi(z_2,y,a)\!\!\sum_{z_1=x+y-z_2}^N\!\!\!\pi(x-z_2,z_1-y,c)\pi(z_1-y,N-y,b).
  \label{Kconv5}
\end{equation}
For the following $\bm{\lambda}$,
\begin{equation}
  \bm{\lambda}=p\eqdef\frac{bc}{1-a+bc},
\end{equation}
the symmetry condition \eqref{Kcond1} is satisfied and
{\bf Theorem\,\ref{theo1}} gives \eqref{Keig}.
By writing down the eigenvalue formula \eqref{eigform}, in which only the
$z_2=0$ term contributes, we have
\begin{equation}
  \kappa(n)=\sum_{z=0}^N\pi(z,N,b)\sum_{y=0}^z\pi(y,z,c)\check{P}_n(y,p)
  =(a-bc)^n
  ={}_1F_0\Bigl(\genfrac{}{}{0pt}{}{-n}{-}\Bigm|bcp^{-1}\Bigr),
\end{equation}
by using the generating function \eqref{KLS(9.11.11)} and the self-duality (K).
It is interesting to note that these expressions are all symmetric in $b$ and
$c$. Odd eigenvalues are negative if $a<bc$.

\paragraph{Krawtchouk\,$\to$\,Charlier}

This is achieved by fixing $a$ and $c$ with $b\to bN^{-1}$, $N\to\infty$,
\begin{align}
  &\check{P}_n(x,p)\to\check{P}_{\text{C}\,n}(x,p'),\quad
  p'\eqdef\frac{bc}{1-a},\n
  &\pi(x,N,p)\to\pi_{\text{C}}(x,p'),\quad
  \kappa(n)\to\kappa_{\text{C}}(n)=a^n,\n
  &K(x,y)\to K_{\text{C}}(x,y,a,b,c)=\!\!\!\sum_{z_2=0}^{\min(x,y)}\!\!\!
  \pi(z_2,y,a)\!\!\!\sum_{z_1=x+y-z_2}^{\infty}\!\!\!\!
  \pi(x-z_2,z_1-y,c)\pi_{\text{C}}(z_1-y,b),
  \label{KCconv5}
\end{align}
and the relations \eqref{Kcond1}, \eqref{Keig} and \eqref{eigform} of (K)
reduce to those of (C).
Based on $K_{\text{C}}$ \eqref{KCconv5}, let us rederive these results.
The symmetry condition \eqref{Kcond1inf} is satisfied for $\bm{\lambda}=p'$ and
{\bf Theorem\,\ref{theo1}} gives \eqref{Keiginf}.
The eigenvalue formula \eqref{eigforminf} reads
\begin{equation}
  \kappa_{\text{C}}(n)
  =\sum_{z=0}^{\infty}\pi_{\text{C}}(z,b)\sum_{y=0}^z\pi(y,z,c)
  \check{P}_{\text{C}\,n}(y,p')=a^n
  ={}_1F_0\Bigl(\genfrac{}{}{0pt}{}{-n}{-}\Bigm|bcp^{\prime\,-1}\Bigr),
\end{equation}
by using the generating function \cite{kls}(9.14.11) and the self-duality (C).

\subsection{Multiple convolutions of type (\romannumeral1) and (\romannumeral3)}
\label{sec:(i)(iii)gen}

For the Krawtchouk (K) case, type (\romannumeral1) and (\romannumeral3)
convolutions can be repeated indefinitely.

Let us define $\pi^{(\pm)}(x,y,N,p)$ ($x,y\in\mathbb{Z}$) as follows:
\begin{align*}
  \pi^{(+)}(x,y,N,p)&\eqdef
  \left\{\begin{array}{ll}
  \pi(x-y,N-y,p)&:0\le y\le x\le N\\
  0&:\text{otherwise}
  \end{array}\right.,\\
  \pi^{(-)}(x,y,N,p)&\eqdef
  \left\{\begin{array}{ll}
  \pi(x,y,p)&:0\le x\le y\le N\\
  0&:\text{otherwise}
  \end{array}\right..
\end{align*}
They are related by
\begin{equation}
  \pi^{(-)}(x,y,N,p)=\pi^{(+)}(N-x,N-y,N,1-p).
  \label{pi+-}
\end{equation}
For an integer $m\ge2$, let us define
$K(x,y)=K^{(\epsilon_1,\ldots,\epsilon_m)}(x,y,N,p_1,\ldots,p_m)$
($x,y\in\cX$, $\epsilon_j=\pm$) by
\begin{equation}
  K(x,y)\eqdef\sum_{z_1,\ldots,z_{m-1}=0}^N\,\prod_{j=1}^m
  \pi^{(\epsilon_j)}(z_{j-1},z_j,N,p_j)
  \quad(z_0=x,\ z_m=y).
  \label{Kxym}
\end{equation}
For example,
\begin{align*}
  &\quad K^{(+,-,+,-)}(x,y,N,p_1,p_2,p_3,p_4)\n
  &=\sum_{z_2=0}^{N}\sum_{z_1=0}^{\min(x,z_2)}\sum_{z_3=0}^{\min(z_2,y)}
  \pi(x-z_1,N-z_1,p_1)\pi(z_1,z_2,p_2)\pi(z_2-z_3,N-z_3,p_3)\pi(z_3,y,p_4).
\end{align*}
For $m=2$ case, $K^{(+,-)}(x,y,N,p_1,p_2)$ and $K^{(-,+)}(x,y,N,p_1,p_2)$
correspond to \eqref{KK1} and \eqref{KK3} with $(a,b)=(p_2,p_1)$, respectively.
{}From \eqref{pi+-}, we have
\begin{equation}
  K^{(\epsilon_1,\ldots,\epsilon_m)}(x,y,N,p_1,\ldots,p_m)
  =K^{(-\epsilon_1,\ldots,-\epsilon_m)}(N-x,N-y,N,1-p_1,\ldots,1-p_m),
  \label{Kep=K-e1-p}
\end{equation}
and \eqref{KK1} and \eqref{KK3} are connected by this relation.
Since two successive $\pi^{(+)}\pi^{(+)}$ and $\pi^{(-)}\pi^{(-)}$ can be
reduced to one $\pi^{(+)}$ and $\pi^{(-)}$ by \eqref{Kpipi(ii)} and
\eqref{Kpipi(i)}, respectively, it is sufficient to consider
$(\epsilon_1,\ldots,\epsilon_m)=(+,-,+,-,\ldots)$ or $(-,+,-,+,\ldots)$.
That is, the multiple convolutions of type (\romannumeral1) and (\romannumeral3),
respectively.
Since $\pi^{(+)}$ is lower triangular and $\pi^{(-)}$ is upper triangular,
\begin{align*}
  \det\pi^{(+)}(*,*,N,p)&=\prod_{x\in\cX}\pi^{(+)}(x,x,N,p)
  =\prod_{n\in\cX}(1-p)^n,\\
  \det\pi^{(-)}(*,*,N,p)&=\prod_{x\in\cX}\pi^{(-)}(x,x,N,p)=\prod_{n\in\cX}p^n,
\end{align*}
the eigenvalues are easily guessed as in the original type (\romannumeral1)
and (\romannumeral3) cases,
\begin{equation}
  \det K^{(\epsilon_1,\ldots,\epsilon_m)}
  =\prod_{n\in\cX}\biggl(\prod_{j=1}^mp_j^{(\epsilon_j)}\biggr)^n,\quad
  p^{(+)}\eqdef 1-p,\quad p^{(-)}\eqdef p.
\end{equation}
In order to determine the parameter $\bm{\lambda}=p$ which satisfies the
symmetry condition \eqref{Kcond1}, we solve the equation
\begin{equation*}
  K(0,N)\pi(N,N,p)=K(N,0)\pi(0,N,p)\quad\Biggl(\Rightarrow
  p=\frac{1}{1+\Bigl(\frac{K(0,N)}{K(N,0)}\Bigr)^{\frac{1}{N}}}\Biggr),
\end{equation*}
which is obtained by setting $x=0$ and $y=N$ in \eqref{Kcond1}.
Among $(N+1)\times(N+1)$ elements of $K(x,y)$, $K(0,N)$ and $K(N,0)$ are the
easiest to evaluate, by successive applications of the binomial theorem,
for example,
\begin{align*}
  &K^{(+,-)}(0,N,N,p_1,p_2)=(1-p_1)^N(1-p_2)^N,\quad
  K^{(+,-)}(N,0,N,p_1,p_2)=p_1^N,\\
  &K^{(+,-,+)}(0,N,N,p_1,p_2,p_3)=(1-p_1)^N(1-p_2)^N,\\ 
  &K^{(+,-,+)}(N,0,N,p_1,p_2,p_3)=\bigl(p_1+(1-p_1)p_2p_3\bigr)^N,\\
  &K^{(+,-,+,-)}(0,N,N,p_1,p_2,p_3,p_4)
  =\bigl((1-p_1)(1-p_2p_3-p_2p_4+p_2p_3p_4)\bigr)^N,\\
  &K^{(+,-,+,-)}(N,0,N,p_1,p_2,p_3,p_4)=\bigl(p_1+(1-p_1)p_2p_3\bigr)^N.
\end{align*}
We obtain the following $\bm{\lambda}=p$ and
\begin{align}
  \kappa(n)&=\biggl(\prod_{j=1}^mp_j^{(\epsilon_j)}\biggr)^n
  ={}_1F_0\Bigl(\genfrac{}{}{0pt}{}{-n}{-}\Bigm|1-\kappa(1)\Bigr),\\
  (\epsilon_1,\ldots,\epsilon_m)&=(+,-,+,-,\ldots):
  \ \ p=\frac{1}{1-\kappa(1)}\sum_{k=0}^{[\frac{m-1}{2}]}
  \prod_{j=1}^{2k}p_j^{(\epsilon_j)}\cdot p_{2k+1},\\
  (\epsilon_1,\ldots,\epsilon_m)&=(-,+,-,+,\ldots):
  \ \ p=\frac{1}{1-\kappa(1)}\sum_{k=1}^{[\frac{m}{2}]}
  \prod_{j=1}^{2k-1}p_j^{(\epsilon_j)}\cdot p_{2k}.
\end{align}
These two are related by \eqref{Kep=K-e1-p}.
The symmetry condition \eqref{Kcond1} is verified by explicit calculation
for small $m$ and $N$. We do not have an analytical proof of the symmetry
condition \eqref{Kcond1} for general $m$ and $N$.

\subsection{One-parameter families of commuting $K$'s}
\label{sec:powerL}

We have derived various $K(x,y)$'s satisfying
\begin{align*}
  &K(x,y)>0,\quad\sum_{x\in\cX}K(x,y)=1,
  \tag{\ref{basK}}\\
  &\sum_{y\in\cX}K(x,y)\pi(y,N,\bm{\lambda})\check{P}_n(y,\bm{\lambda})
  =\kappa(n)\pi(x,N,\bm{\lambda})\check{P}_n(x,\bm{\lambda})\ \ (n\in\cX).
  \tag{\ref{Keig}}
\end{align*}
It is trivial that the $m$-th power of $K$ ($m\ge1$), $K^m$, also satisfies
\begin{align*}
  &K^m(x,y)>0,\quad\sum_{x\in\cX}K^m(x,y)=1,\n
  &\sum_{y\in\cX}K^m(x,y)\pi(y,N,\bm{\lambda})\check{P}_n(y,\bm{\lambda})
  =\kappa(n)^m\pi(x,N,\bm{\lambda})\check{P}_n(x,\bm{\lambda})\ \ (n\in\cX).
\end{align*}
Namely $K^m$ also gives an exactly solvable Markov chain.

It is interesting to note that by changing the parameters $\{\bm{\lambda}_j\}$
to $\{\bm{\lambda}_j'(t)\}$, some $K$'s derived in
\S\,\ref{convtype1}--\S\,\ref{convtype5} can be deformed to create a one
parameter ($t$) family of commuting $K$'s.
That is, they share the same eigenvectors but the eigenvalues are different.
For examples, the following $\bm{\lambda}'_j(t)$'s give the same $\bm{\lambda}$
in the symmetry condition \eqref{Kcond1},
\begin{align*}
  \eqref{K1la}:&\ \ \bm{\lambda}'_1(t)=at,
  \ \ \bm{\lambda}'_2(t)=\frac{(1-at)b}{1-a(1-b(1-t))}\ \ (0<t\le1)
  \Rightarrow\bm{\lambda}=\frac{b}{1-a+ab},\\
  \eqref{KCconv1}:&\ \ \bm{\lambda}'_1(t)=at,
  \ \ \bm{\lambda}'_2(t)=\frac{(1-at)b}{1-a}\ \ (0<t\le1)
  \Rightarrow\bm{\lambda}=\frac{b}{1-a},\\
  \eqref{H1la}:&\ \ \bm{\lambda}'_1(t)=(a+t,b-t),
  \ \ \bm{\lambda}'_2(t)=(b-t,c)\ \ (-a<t<b)
  \Rightarrow\bm{\lambda}=(a+b,c),\\
  \eqref{qH1la}:&\ \ \bm{\lambda}'_1(t)=(at,bt^{-1}),
  \ \ \bm{\lambda}'_2(t)=(bt^{-1},c)\ \ (b<t<a^{-1})
  \Rightarrow\bm{\lambda}=(ab,c),\\
  \eqref{KqMr1}:&\ \ \bm{\lambda}'_1(t)=(at,bt^{-1}),
  \ \ \bm{\lambda}'_2(t)=(bt^{-1}q^{-1},b^{-1}tc)\ \ (b<t<a^{-1})
  \Rightarrow\bm{\lambda}=\bigl(abq^{-1},(ab)^{-1}c\bigr).
\end{align*}
The two matrices $K(x,y,\{\bm{\lambda}_j\})$ and $K(x,y,\{\bm{\lambda}'_j(t)\})$
have the common eigenvectors $\pi(x,\bm{\lambda})$ $\check{P}_n(x,\bm{\lambda})$
and they commute as is clear from the spectral representation
{\bf Theorem\,\ref{theo3}}.
But the eigenvalues are different,
$\kappa(n,\{\bm{\lambda}_j\})\neq\kappa(n,\{\bm{\lambda}'_j(t)\})$.
Let $\bm{\lambda}_j'(t^{[i]})$ ($i=1,2,\ldots$) be the parameters giving the
same $\bm{\lambda}$, and set $K^{[i]}(x,y)=K(x,y,\{\bm{\lambda}_j'(t^{[i]})\})$
and $\kappa^{[i]}(n)=\kappa(n,\{\bm{\lambda}_j'(t^{[i]})\})$.
For $m$ such $K^{[i]}$'s, let us consider their matrix product
(the order is irrelevant),
\begin{equation}
  K^{(m)}\eqdef K^{[m]}\cdots K^{[2]}K^{[1]},\quad
  \kappa^{(m)}(n)\eqdef\kappa^{[m]}(n)\cdots\kappa^{[2]}(n)\kappa^{[1]}(n).
\end{equation}
Then we have
\begin{align}
  &K^{(m)}(x,y)>0,\quad\sum_{x\in\cX}K^{(m)}(x,y)=1,\n
  &\sum_{y\in\cX}K^{(m)}(x,y)\pi(y,N,\bm{\lambda})\check{P}_n(y,\bm{\lambda})
  =\kappa^{(m)}(n)\pi(x,N,\bm{\lambda})\check{P}_n(x,\bm{\lambda})\ \ (n\in\cX).
\end{align}
Namely $K^{(m)}$ also gives an exactly solvable Markov chain.

\section{Other Topics}
\label{sec:other}

Here we discuss two related topics.

\subsection{Dual Markov chains}
\label{sec:dual}

For a Markov chain $K(x,y)$ on a finite one dimensional integer lattice $\cX$,
its `dual' Markov chain $K^{\text{d}}(x,y)$ is defined by the similarity
transformation in terms of the anti-diagonal matrix $J$,
$J(x,y)\eqdef\delta_{x,N-y}$,
\begin{align} 
  &\!K^{\text{d}}(x,y)\eqdef(JKJ)(x,y)=K(N-x,N-y),\\
  &\!\sum_{y\in\cX}K(x,y)v_n(y)=\kappa(n)v_n(x)\,\Rightarrow
  \sum_{y\in\cX}K^{\text{d}}(x,y)v^{\text{d}}_n(y)=\kappa(n)v^{\text{d}}_n(x),
  \ v^{\text{d}}_n(x)\eqdef v_n(N-x).\!\!
\end{align}
The concept of duality was reported in \cite{coo-hoa-rah77}.

For $K(x,y)$ constructed by convolutions listed in \S\,\ref{sec:convlist},
the dual Markov chains take the following forms,
\begin{align} 
  \text{(d\romannumeral1)}:&\ \ K^{\text{d}}(x,y)\eqdef\sum_{z=\max(x,y)}^N
  \!\!\!\pi(z-x,z,\bm{\lambda}_2)\pi(N-z,N-y,\bm{\lambda}_1),
  \label{dconv11}\\
  \text{(d\romannumeral2)}:&\ \ K^{\text{d}}(x,y)\eqdef
  \sum_{z=\max(x,y)}^{\min(x+y,N)}
  \!\!\!\pi(z-x,y,\bm{\lambda}_2)\pi(N-z,N-y,\bm{\lambda}_1),
  \label{dconv2}\\
  \text{(d\romannumeral3)}:&\ \ K^{\text{d}}(x,y)\eqdef\sum_{z=0}^{\min(x,y)}
  \!\!\!\pi(N-x,N-z,\bm{\lambda}_2)\pi(y-z,y,\bm{\lambda}_1),
  \label{dconv3}\\
  \text{(d\romannumeral4)}:&\ \ K^{\text{d}}(x,y)\eqdef
  \!\!\!\!\!\sum_{z_2=\max(x,y)}^N\!\!\!\!\!\!
  \pi(N-z_2,N-y,\bm{\lambda}_1)\!\!\!\sum_{z_1=0}^{\min(x,y)}\!\!\!
  \pi(z_2-x,z_2-z_1,\bm{\lambda}_3)\pi(y-z_1,y,\bm{\lambda}_2),
  \label{dconv4}\\
  \text{(d\romannumeral5)}:&\ \ K^{\text{d}}(x,y)\eqdef
  \!\!\!\!\!\sum_{z_2=\max(x,y)}^N\!\!\!\!\!
  \pi(N-z_2,N-y,\bm{\lambda}_1)\!\!\!\sum_{z_1=0}^{x+y-z_2}\!\!\!
  \pi(z_2-x,y-z_1,\bm{\lambda}_3)\pi(y-z_1,y,\bm{\lambda}_2),
  \label{dconv5}\\
  &\!\!\!\!\!\!\!\!\!\!\!\! K^{\text{d}}(x,y)\pi(N-y,N,\bm{\lambda})
  =K^{\text{d}}(y,x)\pi(N-x,N,\bm{\lambda}).
\end{align}
For (K) and (H), the dual eigenvectors take the standard forms with flipped
$\bm{\lambda}$, see \eqref{Kpid}, \eqref{Kpd} for (K) and \eqref{Hpi},
\eqref{Hpd} for (H),
\begin{alignat*}{2}
  \pi_{\text{K}}(N-x,N,p)&=\pi_{\text{K}}(x,N,1-p),
  &\ \ \check{P}_{\text{K}\,n}(N-x,p)&\propto\check{P}_{\text{K}\,n}(x,1-p),\\
  \pi_{\text{H}}(N-x,N,a,b)&=\pi_{\text{H}}(x,N,b,a),
  &\ \ \check{P}_{\text{H}\,n}(N-x,a,b)&\propto\check{P}_{\text{H}\,n}(x,b,a),
\end{alignat*}
and the limiting procedure for $N\to\infty$ with fixed $x$ goes in a similar
way as before. However, for ($q$H) the situation is different, see \eqref{qHrn},
\begin{equation}
  \pi_{\text{$q$H}}(N-x,N,a,b)\not\propto\pi_{\text{$q$H}}(x,N,b,a),\quad 
  \check{P}_{\text{$q$H}\,n}(N-x,N,a,b)\not\propto
  \check{P}_{\text{$q$H}\,n}(x,N,b,a),
  \label{piqHneqpiqH}
\end{equation}
and the reductions to ($q$M) cannot be carried out.

\subsubsection{reduction to little $q$-Jacobi}
\label{sec:lqJ}

It is known \cite{os12} (\S\,V.C.1) that the $q$-Hahn polynomial with the
replacement $x\to N-x$ is another polynomial in $\eta(x)=1-q^x$ called the
alternative $q$-Hahn (a$q$H) polynomial.
The basic data of (a$q$H) with $\bm{\lambda}=(a,b)$ ($0<a<1$, $b<1$) are
\begin{align}
  &\pi(x,N,a,b)=\genfrac{[}{]}{0pt}{}{\,N\,}{x}
  \frac{(a\,;q)_{N-x}\,(b\,;q)_x\,a^x}{(ab\,;q)_N},\\
  &d_n^2=\genfrac{[}{]}{0pt}{}{\,N\,}{n}
  \frac{(b,abq^{-1}\,;q)_n\,a^nq^{n(n-1)}}{(a,abq^N\,;q)_n}
  \frac{1-abq^{2n-1}}{1-abq^{-1}},\\
  &s_1\eta(x)\pi(x,N,a,b)=-\pi(x-1,N-1,a,bq),\quad
  s_1\eqdef-\frac{1-ab}{a(1-b)(1-q^N)},\\
  &\eta(z)\pi(z,x,a,b)=\frac{a(1-b)}{1-ab}\eta(x)\pi(z-1,x-1,a,bq),\quad
  \eta(x)=1-q^x,\\
  &\check{P}_n(x,a,b)=P_n\bigl(\eta(x),a,b\bigr)
  ={}_3\phi_2\Bigl(\genfrac{}{}{0pt}{}{q^{-n},\,abq^{n-1},\,q^{-x}}
  {b,\,q^{-N}}\Bigm|q\,;a^{-1}q^{x+1-N}\Bigr).
\end{align}
The above non-proportionality relation \eqref{piqHneqpiqH} is rewritten as the
equalities between the $q$-Hahn ($q$H) and the alternative $q$-Hahn (a$q$H)
polynomials
\begin{align}
  \pi_{\text{$q$H}}(N-x,N,a,b)&=\pi_{\text{a$q$H}}(x,N,a,b),\\
  \check{P}_{\text{$q$H}\,n}(N-x,a,b)
  &=(-a)^nq^{\binom{n}{2}}\frac{(b\,;q)_n}{(a\,;q)_n}
  \check{P}_{\text{a$q$H}\,n}(x,a,b).
\end{align}
The transition matrices $K$ of (a$q$H) are expressed by the duals of those of
($q$H), see \eqref{dconv11}--\eqref{dconv5},
\begin{align}
  K^{\text{(\romannumeral1)}}_{\text{a$q$H}}
  (x,y,\bm{\lambda}_1,\bm{\lambda}_2)
  &=K^{\text{(d\romannumeral3)}}_{\text{$q$H}}
  (x,y,\bm{\lambda}_1,\bm{\lambda}_2),
  \label{KiaqH=KiiiqH}\\
  K^{\text{(\romannumeral3)}}_{\text{a$q$H}}
  (x,y,\bm{\lambda}_1,\bm{\lambda}_2)
  &=K^{\text{(d\romannumeral1)}}_{\text{$q$H}}
  (x,y,\bm{\lambda}_1,\bm{\lambda}_2),
  \label{KiiiaqH=KiqH}\\
  K^{\text{(\romannumeral4)}}_{\text{a$q$H}}
  (x,y,\bm{\lambda}_1,\bm{\lambda}_2,\bm{\lambda}_3)
  &=K^{\text{(d\romannumeral4)}}_{\text{$q$H}}
  (x,y,\bm{\lambda}_2,\bm{\lambda}_1,\bm{\lambda}_3),
  \label{KivaqH=KivqH}
\end{align}
independently of the choice of parameters $\bm{\lambda}_1$ and $\bm{\lambda}_2$
(and $\bm{\lambda}_3$).
Therefore the results of ($q$H) in \S\,\ref{1qH}, \S\,\ref{3qH} and \S\,\ref{4qH}
are translated to those of (a$q$H).
This is merely rewriting, but their $N\to\infty$ limits give new results.
By the limit $N\to\infty$ with fixed $a$ and $b$, (a$q$H) goes to the little
$q$-Jacobi (l$q$J) polynomial,
\begin{equation}
  \lim_{N\to\infty}\pi_{\text{a$q$H}}(x,N,a,b)
  =\pi_{\text{l$q$J}}(x,a,b),\quad
  \lim_{N\to\infty}\check{P}_{\text{a$q$H}\,n}(x,N,a,b)
  =\check{P}_{\text{l$q$J}\,n}(x,a,b).
\end{equation}
The basic data of (l$q$J) with $\bm{\lambda}=(a,b)$ ($0<a<1$, $b<1$) is
\begin{align}
  &\pi(x,a,b)=\frac{(b\,;q)_x\,a^x}{(q\,;q)_x}
  \frac{(a\,;q)_{\infty}}{(ab\,;q)_{\infty}},\quad
  d_n^2=\frac{(b,abq^{-1}\,;q)_n\,a^nq^{n(n-1)}}{(q,a\,;q)_n}
  \frac{1-abq^{2n-1}}{1-abq^{-1}},
  \label{pilqJ}\\
  &s_1\eta(x)\pi(x,a,b)=-\pi(x-1,a,bq),\quad\eta(x)=1-q^x,\quad
  s_1\eqdef-\frac{1-abq^{-1}}{a(1-b)},
  \label{lqJtrif1}\\
  &\check{P}_n(x,a,b)=P_n\bigl(\eta(x),a,b\bigr)
  ={}_3\phi_1\Bigl(\genfrac{}{}{0pt}{}{q^{-n},\,abq^{n-1},\,q^{-x}}{b}
  \Bigm|q\,;a^{-1}q^{x+1}\Bigr)\\
  &\phantom{\check{P}_n(x,a,b)=P_n\bigl(\eta(x),a,b\bigr)}
  =(-a)^{-n}q^{-\binom{n}{2}}\frac{(a\,;q)_n}{(b\,;q)_n}\,
  {}_2\phi_1\Bigl(\genfrac{}{}{0pt}{}{q^{-n},\,abq^{n-1}}{a}
  \Bigm|q\,;q^{x+1}\Bigr).
\end{align}
The conventional little $q$-Jacobi polynomial is
$p_n(q^x;a,b|q)={}_2\phi_1\bigl(\genfrac{}{}{0pt}{}{q^{-n},\,abq^{n+1}}{aq}
\bigl|q\,;q^{x+1}\bigr)$
and our parametrisation is slightly different from the standard one
$(a,b)^{\text{standard}}=(aq^{-1},bq^{-1})$.
Similarly to those examples in \S\ref{sec:selfsim}, this $\pi$ keeps its form
under the following convolutions:
\begin{align}
  &\sum_{z=0}^x\pi(x-z,a_1,b_1)\pi(z,a_1b_1,b_2)=\pi(x,a_1,b_1b_2),\\
  &\sum_{z=y}^x\pi(x-z,a_1,b_1)\pi(z-y,a_1b_1,b_2)=\pi(x-y,a_1,b_1b_2).
\end{align}
These are obtained by the sum formula for $\pi_{\text{$q$H}}$ \eqref{bt4}.

\paragraph{alternative $q$-Hahn\,$\to$\,little $q$-Jacobi}

{}From the type (\romannumeral1), (\romannumeral3), (\romannumeral4)
transition matrices $K_{\text{a$q$H}}$, we obtain $K_{\text{l$q$J}}$
whose eigenvectors are described by the little $q$-Jacobi polynomial,
\begin{align}
  \text{(\romannumeral1)}
  &:\ K_{\text{l$q$J}}(x,y)=\sum_{z=0}^{\min(x,y)}\!\!
  \pi_{\text{l$q$J}}(x-z,c,a)\pi_{\text{a$q$H}}(z,y,a,b),
  \label{(i)lqJ}\\
  &\quad\ \text{eigenvector}:\pi_{\text{l$q$J}}(x,c,ab)
  \check{P}_{\text{l$q$J}\,n}(x,c,ab),\quad
  \kappa(n):\eqref{3qHphi},\\
  \text{(\romannumeral3)}
  &:\ K_{\text{l$q$J}}(x,y)=\!\!\sum_{z=\max(x,y)}^{\infty}\!\!\!\!
  \pi_{\text{a$q$H}}(x,z,b,c)\pi_{\text{l$q$J}}(z-y,a,b),
  \label{(iii)lqJ}\\
  &\quad\ \text{eigenvector}:\pi_{\text{l$q$J}}(x,ab,c)
  \check{P}_{\text{l$q$J}\,n}(x,ab,c),\quad
  \kappa(n):\eqref{1qHphi},\\
  \text{(\romannumeral4)}
  &:\ K_{\text{l$q$J}}(x,y)=\sum_{z_2=0}^{\min(x,y)}\!\!
  \pi_{\text{a$q$H}}(z_2,y,a_2,b_2)\!\!
  \sum_{z_1=\max(x,y)}^{\infty}\!\!\!\!
  \pi_{\text{a$q$H}}(x-z_2,z_1-z_2,b_1,a_2)\n
  &\hspace{78mm}\times
  \pi_{\text{l$q$J}}(z_1-y,a_1,b_1),
  \label{(iv)lqJ}\\
  &\quad\ \text{eigenvector}:
  \pi_{\text{l$q$J}}(x,a_1b_1,a_2b_2)
  \check{P}_{\text{l$q$J}\,n}(x,a_1b_1,a_2b_2),\quad
  \kappa(n):\eqref{(4.109)}.
  \label{(iv)lqJeigen}
\end{align}

\subsection{Repeated discrete time Birth and Death processes}
\label{sec:repBD}

Exactly solvable discrete time Birth and Death (BD) processes $K_{\text{BD}}$
are constructed \cite{dtbd} based on exactly solvable continuous time BD on a
one dimensional integer lattice $\cX$,
\begin{align}
  &K_{\text{BD}}=I+t_{\text{S}}\,L_{\text{BD}},
  \ \text{i.e.}\ K_{\text{BD}}(x,y)=\delta_{x\,y}
  +t_{\text{S}}\,L_{\text{BD}}(x,y),\\
  &L_{\text{BD}}(x+1,x)=B(x),\ \ L_{\text{BD}}(x-1,x)=D(x),
  \ \ L_{\text{BD}}(x,x)=-B(x)-D(x),\n
  &\quad L_{\text{BD}}(x,y)=0\ \ (|x-y|\ge2),\quad
  \sum_{x\in\cX}L_{\text{BD}}(x,y)=0,
  \label{LBDdef}
\end{align}
in which $B(x)$ and $D(x)$ are the birth and death rates at point $x$ and they
are chosen to be the coefficient functions of the difference equations governing
the orthogonal polynomials $\{\check{P}_n(x)\}$ belonging to Askey scheme
\cite{bdsol,askey,kls,os12},
\begin{align}
  B(x)\bigl(\check{P}_n(x)-\check{P}_n(x+1)\bigr)
  +D(x)\bigl(\check{P}_n(x)-\check{P}_n(x-1)\bigr)
  =\mathcal{E}(n)\check{P}_n(x)\quad(n\in\cX),
  \label{bdeq0}
\end{align}
and the time scale parameter $t_{\text{S}}$ must satisfy the upper bound
condition
\begin{equation}
  t_{\text{S}}\cdot\max\bigl(B(x)+D(x)\bigr)<1.
  \label{dctrel}
\end{equation}
The eigenvectors of $K_{\text{BD}}$ are $\{\pi(x)\check{P}_n(x)\}$ ($n\in\cX$)
and $\pi(x)$ is the normalised orthogonality measure of the polynomial $=$
the stationary probability distribution and
\begin{equation}
  \sum_{y\in\cX}K_{\text{BD}}(x,y)\pi(y)\check{P}_n(y)
  =\kappa(n)\pi(x)\check{P}_n(x),\quad
  \kappa(n)\eqdef 1-t_{\text{S}}\,\mathcal{E}(n)\quad(n\in\cX).
\end{equation}
This applies to a good part of the orthogonal polynomials of a discrete variable
in Askey scheme, including the Hahn, $q$-Hahn and Racah and $q$-Racah
\cite{dtbd}.

{}From the definition \eqref{LBDdef}, it is shown that the $m$-th power of
$L_{\text{BD}}$ ($m\ge1$), $L_{\text{BD}}^m$, has the following form of
the matrix elements,
\begin{equation*}
  L_{\text{BD}}^m(x+k,x)=(-1)^{m-k}a^{(m)}_k(x)\ \ (-m\le k\le m),\quad
  L_{\text{BD}}^m(x,y)=0\ \ (|x-y|>m),
\end{equation*}
where $a^{(m)}_k(x)>0$.
Let us consider the following matrix $X$,
\begin{equation}
  X=\sum_{j=0}^{m-1}c_jL_{\text{BD}}^{m-j},\quad c_0=1\quad
  \Bigl(\Rightarrow\sum_{x\in\cX}X(x,y)=0,\ X(x,y)=0\ (|x-y|>m)\Bigr),
  \label{Xdef}
\end{equation}
where $c_j$ are constants. Its non zero matrix elements are
\begin{align*}
  X\bigl(x\pm(m-k),x\bigr)&=\sum_{j=0}^kc_j(-1)^{k-j}a^{(m-j)}_{\pm(m-k)}(x)
  \ \ (0\le k\le m-1),\\
  X(x,x)&=\sum_{j=0}^{m-1}c_j(-1)^{m-j}a^{(m-j)}_0(x).
\end{align*}
Starting from $X(x\pm m,x)=a^{(m)}_{\pm m}(x)>0$, we can tune $c_k$
($k=1,\ldots,m-2$ in turn) such that $X(x\pm(m-k),x)>0$, and tune $c_{m-1}$
such that $X(x\pm 1,x)>0$ and $X(x,x)<0$.
For such chosen weights $\{c_j\}$, we define a matrix $K^{(m)}_{\text{BD}}$,
\begin{equation}
  K^{(m)}_{\text{BD}}\eqdef I_d+t_{\text{S}}\,X,\quad
  t_{\text{S}}\cdot\max\bigl(-X(x,x)\bigr)<1,
\end{equation}
which satisfies
\begin{equation*}
  K^{(m)}_{\text{BD}}(x,y)\ge0,\quad
  \sum_{x\in\cX}K^{(m)}_{\text{BD}}(x,y)=1,\quad
  K^{(m)}_{\text{BD}}(x,y)=0\ \ (|x-y|>m).
\end{equation*}
This gives an exactly solvable Markov chain and
the matrices $K^{(m)}_{\text{BD}}$'s have common eigenvectors
\begin{align}
  &\sum_{y\in\cX}K^{(m)}_{\text{BD}}(x,y)\pi(y)\check{P}_n(y)
  =\kappa^{(m)}(n)\pi(x)\check{P}_n(x),\n
  &\kappa^{(m)}(n)\eqdef 1+t_{\text{S}}
  \sum_{j=0}^{m-1}(-1)^{m-j}c_j\mathcal{E}(n)^{m-j}\quad(n\in\cX).
\end{align}

\bigskip

Before closing this section a few remarks on possible applications are in order.
As is well known the birth and death processes are recognised to be related to
diffusion processes \cite{feller}. In other words, the equations for BD
processes are space discretisation of 1-d Fokker-Planck equations.
Discrete time BD processes are further discretisation in time.
As illustrated by the famous Ehrenfest urn model, Markov chains are also
related to diffusion processes as well as to the familiar random walks.
Various examples in this and preceding sections are multi-parameter
generalisations of known Markov chains and BD processes.
It is expected that they would find diverse applications in physics,
chemistry, etc. in particular, diffusions and random walks.

\section{Summary and Comments}
\label{sec:comm}

Based on the fact that many examples of exactly solvable BD processes/chains
\cite{bdsol,dtbd}, typical cases of Markov processes/chains, have been
constructed in terms of orthogonal polynomials of a discrete variable
\cite{nikiforov}--\cite{gasper}, we establish a wide range of generalisations
of various Markov chains/processes \cite{hoa-rah83}--\cite{diaconis20}
using solvability as a guide. Adopting the convolutions of the orthogonality
measures is a new input of the present research.
It is an interesting challenge to formulate
probabilistic procedures/interpretations like ``cumulative Bernoulli trials''
for these new examples.

Except for those corresponding to the extra eigenvectors for the $q$-Meixner
($q$M) in \S\,\ref{1qH}, \S\,\ref{3qH}, \S\,\ref{4qH}, the eigenvalues of the
$K$'s derived in this paper share many remarkable properties;
\begin{enumerate}
\item
$\kappa(n)$ has a neat sum formula of one or two $\pi$'s and $\check{P}_n$,
for example
\begin{equation*}
  \kappa(n)=\sum_{y\in\cX}\pi(y,N,\bm{\lambda}')\check{P}_n(y,\bm{\lambda}).
\end{equation*}
\item
$\kappa(n)$ is independent of $N$, the size of the lattice.
This reminds us of the similar situation of the exactly solvable BD processes
\cite{bdsol,os34,dtbd}, which is due to the fact that the eigenvalues of the
difference equations governing these polynomials are $N$ independent
\cite{os12}.
\item
$\kappa(n)$ has an expression of a terminating ($q$)-hypergeometric series,
for example,
\begin{equation*}
  \kappa(n)={}_3\phi_2\Bigl(\genfrac{}{}{0pt}{}
  {q^{-n},\,b_1,\,a_2}{a_1b_1,\,b_1a_2}\Bigm|q\,;a_1b_1q^n\Bigr).
\end{equation*}
We do not know how the bounds $-1<\kappa(n)\le1$ are ingrained in the
hypergeometric expressions, in particular, in ${}_3F_2$, ${}_4F_3$,
${}_3\phi_2$ and ${}_4\phi_3$.
\end{enumerate}
Deciphering these curious facts, we believe, would lead to a deeper
understanding of the subject.

\section*{Acknowledgements}

S.\,O. is supported by JSPS KAKENHI Grant Number JP19K03667.

\bigskip
\section*{Appendix}
\label{sec:app}
\renewcommand{\theequation}{A.\arabic{equation}}

Here we provide the proof of triangularity {\bf Lemma}
\begin{align*}
  \sum_{y\in\cX}K(y,x)\eta(y)^n&=\sum_{m=0}^na_{n\,m}\eta(x)^m\ \ (n\in\cX)
  \quad\bigl(a_{n\,m}=0\ \,\text{for}\ \,n<m\bigr),
  \tag{\ref{triang1}}\\
  \eta(x)&=\left\{
  \begin{array}{ll}
  x&:\text{(\romannumeral1)--(\romannumeral5)}\\
  q^{-x}-1&:\text{(\romannumeral1)},\,\text{(\romannumeral3)},
  \,\text{(\romannumeral4)}
  \end{array}\right..
  \tag{\ref{triang2}}
\end{align*}
The proof depends on a universal property of the normalised orthogonality
measure $\pi(x,N,\bm{\lambda})$ (or $\pi(x,\bm{\lambda})$) of all the
polynomials of a discrete variable in Askey scheme \cite{os12},
{\em except for those having the Jackson integral measures}:
\begin{alignat}{2}
  &\text{finite}:&s_1\eta(x,\bm{\lambda})\pi(x,N,\bm{\lambda})
  &=-\pi(x-1,N-1,\bm{\lambda}'),
  \label{pieta1}\\
  &\text{semi-infinite}:&s_1\eta(x,\bm{\lambda})\pi(x,\bm{\lambda})
  &=-\pi(x-1,\bm{\lambda}'),
  \label{pieta2}
\end{alignat}
in which $s_1$ is the coefficient of $\eta(x,\bm{\lambda})$ in
$P_1(\eta(x,\bm{\lambda}),\bm{\lambda})$,
\begin{equation*}
  P_1\bigl(\eta(x,\bm{\lambda}),\bm{\lambda}\bigr)=1+s_1\eta(x,\bm{\lambda}).
\end{equation*}
Note that the Racah, dual Hahn etc.\,\cite{askey,ismail} are polynomials in
$\eta(x,\bm{\lambda})$, depending on parameters $\bm{\lambda}$.
It is easy to verify these formulas one by one.
In \S\,\ref{sec:data} the explicit expressions of $s_1$ and $\bm{\lambda}'$
are given in \eqref{Ktrif1}, \eqref{Ctrif1}, \eqref{Htrif1}, \eqref{Mtrif1},
\eqref{qHtrif1}, \eqref{qMtrif11}, \eqref{qMtrif21} for (K), (C), (H), (M),
($q$H) and ($q$M), and \eqref{lqJtrif1} for (l$q$J).
Among them, the formulas \eqref{pieta1} for (K), (H) and ($q$H) can be
rewritten as
\begin{equation}
  \eta(z)\pi(z,x,\bm{\lambda})=\beta\,\eta(x)\pi(z-1,x-1,\bm{\lambda}'),
  \label{pieta3}
\end{equation}
in which $\beta$ is a constant independent of $N$.
See \eqref{Ktrif2}, \eqref{Htrif2} and \eqref{qHtrif2} for the explicit forms
for (K), (H) and ($q$H).
To the best of our knowledge, the formulas \eqref{pieta1}--\eqref{pieta3}
have not been reported yet.
 
The general strategy is as follows.
Apply formulas \eqref{pieta1}--\eqref{pieta3} with various arguments,
{\em e.g.}\ $\eta(y-z)$, $\eta(z)$, etc.\ repeatedly to $K(y,x)$ and reduce
$\eta(y)^n$ to a degree $n$ polynomial in $\eta(x)$ through various
intermediaries including $z$ ($z_1$ and $z_2$).
This is guaranteed as the formulas \eqref{pieta1}--\eqref{pieta3} do not
increase the powers of $\eta$.
After the reduction, the remaining $\pi$'s are evaluated to 1 by the summation
in $y$ and $z$ ($z_1$ and $z_2$).
For each type of convolutions $K$, the structure of the reduction
$\eta(y)\to\eta(x)$ is the same for the group having $\eta(x)=x$ {\em i.e.}
(K), (C), (M) and (H). It is more involved for those having $\eta(x)=q^{-x}-1$
{\em i.e.}\ ($q$H) and ($q$M).
It is important to stress that $K$'s for the semi-infinite Markov chains given
in \S\,\ref{sec:exa} have at least one $\pi$ belonging to (K), (H) or ($q$H),
so that \eqref{pieta3} can be applied to extract $\eta(x)$.
The triangularity also holds for these semi-infinite Markov chains.

Below we demonstrate the first step reduction $\eta(y)\to\eta(x)$ for type
(\romannumeral1) convolution
\begin{equation}
  K(y,x)=\sum_{z=0}^{\min(x,y)}
  \!\pi(y-z,N-z,\bm{\lambda}_2)\pi(z,x,\bm{\lambda}_1).
  \label{conv110}
\end{equation}
For $\eta(x)=q^{-x}-1$, by using $\eta(y)=q^{-z}\eta(y-z)+\eta(z)$ and
$q^{-z}\eta(N-z)=\eta(N)-\eta(z)$, we obtain
\begin{align}
  K(y,x)\eta(y)
  &=\sum_{z=0}^{\min(x,y)}\!\!q^{-z}\eta(y-z)\pi(y-z,N-z,\bm{\lambda}_2)
  \pi(z,x,\bm{\lambda}_1)\n
  &\quad+\sum_{z=0}^{\min(x,y)}\!\!\pi(y-z,N-z,\bm{\lambda}_2)
  \eta(z)\pi(z,x,\bm{\lambda}_1)\n
  &=\sum_{z=0}^{\min(x,y-1)}\!\!
  \beta_2q^{-z}\eta(N-z)\pi(y-z-1,N-z-1,\bm{\lambda}'_2)
  \pi(z,x,\bm{\lambda}_1)\n
  &\quad+\sum_{z=1}^{\min(x,y)}\!\!\pi(y-z,N-z,\bm{\lambda}_2)
  \beta_1\eta(x)\pi(z-1,x-1,\bm{\lambda}'_1)\n
  &=\beta_2\eta(N)\!\sum_{z=0}^{\min(x,y-1)}\!\!\!\!
  \pi(y-z-1,N-z-1,\bm{\lambda}'_2)\pi(z,x,\bm{\lambda}_1)\n
  &\quad-\beta_2\beta_1\eta(x)\!\sum_{z=1}^{\min(x,y-1)}\!\!\!\!
  \pi(y-z-1,N-z-1,\bm{\lambda}'_2)\pi(z-1,x-1,\bm{\lambda}'_1)\n
  &\quad+\beta_1\eta(x)\sum_{z=1}^{\min(x,y)}\!\!\pi(y-z,N-z,\bm{\lambda}_2)
  \pi(z-1,x-1,\bm{\lambda}'_1).
 \label{Ketas}
\end{align}
For $\eta(x)=x$, by using $\eta(y)=\eta(y-z)+\eta(z)$ and
$\eta(N-z)=\eta(N)-\eta(z)$, the same result is obtained similarly.
Summing \eqref{Ketas} over $y$, we obtain
\begin{equation*}
  \sum_{y\in\cX}K(y,x)\eta(y)
  =\beta_2\eta(N)+\beta_1(1-\beta_2)\eta(x),
\end{equation*}
which is \eqref{triang1} for $n=1$.

By a similar calculation, we obtain\\
(\romannumeral2) $\eta(x)=x$ :
\begin{align}
  K(y,x)\eta(y)&=\beta_2\eta(N)\!\!\sum_{z=\max(0,x+y-N)}^{\min(x,y-1)}
  \!\!\!\!\!\!\!\!
  \pi(y-z-1,N-x-1,\bm{\lambda}'_2)\pi(z,x,\bm{\lambda}_1)\n
  &\quad+\eta(x)\Bigl(
  -\beta_2\!\!\sum_{z=\max(0,x+y-N)}^{\min(x,y-1)}\!\!\!\!\!\!\!\!
  \pi(y-z-1,N-x-1,\bm{\lambda}'_2)\pi(z,x,\bm{\lambda}_1)\n
  &\phantom{\quad+\eta(x)\Bigl(}
  +\beta_1\!\!\sum_{z=\max(1,x+y-N)}^{\min(x,y)}\!\!\!\!\!\!\!\!
  \pi(y-z,N-x,\bm{\lambda}_2)\pi(z-1,x-1,\bm{\lambda}'_1)\Bigr),
  \label{Keta2}
\end{align}
(\romannumeral3) $\eta(x)=x,q^{-x}-1$ :
\begin{align}
  K(y,x)\eta(y)
  &=\beta_1\beta_2\eta(N)\!\!\sum_{z=\max(x+1,y)}^N\!\!\!\!\!\!
 \pi(y-1,z-1,\bm{\lambda}'_2)\pi(z-x-1,N-x-1,\bm{\lambda}'_1)\n
  &\quad+\beta_2\eta(x)\Bigl(
  -\beta_1\!\sum_{z=\max(x+1,y)}^N\!\!\!\!\!\!
  \pi(y-1,z-1,\bm{\lambda}'_2)\pi(z-x-1,N-x-1,\bm{\lambda}'_1)\n
  &\phantom{\quad+\beta_2\eta(x)\Bigl(}
  +\!\sum_{z=\max(1,x,y)}^N\!\!\!\!\!\!
  \pi(y-1,z-1,\bm{\lambda}'_2)\pi(z-x,N-x,\bm{\lambda}_1)\Bigr),
  \label{Keta3}
\end{align}
(\romannumeral4) $\eta(x)=x,q^{-x}-1$ :
\begin{align}
  &\quad K(y,x)\eta(y)\n
  &=\beta_2\beta_3\eta(N)\!\sum_{z_2=0}^{\min(x,y-1)}\!\!\!\!
  \pi(z_2,x,\bm{\lambda}_1)\!\!\sum_{z_1=\max(x+1,y)}^N\!\!\!\!\!\!\!\!
  \pi(y-z_2-1,z_1-z_2-1,\bm{\lambda}'_3)\n
  &\hspace{66mm}
  \times\pi(z_1-x-1,N-x-1,\bm{\lambda}'_2)\n
  &\quad+\eta(x)\Bigl(
  -\beta_2\beta_3\!\sum_{z_2=0}^{\min(x,y-1)}\!\!\!\!
  \pi(z_2,x,\bm{\lambda}_1)\!\!\sum_{z_1=\max(x+1,y)}^N\!\!\!\!\!\!\!\!
  \pi(y-z_2-\!1,z_1-z_2-\!1,\bm{\lambda}'_3)\n
  &\hspace{74mm}\times
  \pi(z_1-x-\!1,N\!-\!x-\!1,\bm{\lambda}'_2)\n
  &\qquad\quad+\beta_1\sum_{z_2=1}^{\min(x,y)}
  \pi(z_2-1,x-1,\bm{\lambda}'_1)\!\!\sum_{z_1=\max(x,y)}^N\!\!\!\!\!\!
  \pi(y-z_2,z_1-z_2,\bm{\lambda}_3)\pi(z_1-x,N-x,\bm{\lambda}_2)\n
  &\qquad\quad+\beta_3\!\sum_{z_2=0}^{\min(x,y-1)}\!\!\!\!
  \pi(z_2,x,\bm{\lambda}_1)\!\!\!\!\sum_{z_1=\max(x-1,y-1)}^{N-1}
  \!\!\!\!\!\!\!\!\!\!\pi(y-1-z_2,z_1-z_2,\bm{\lambda}'_3)
  \pi(z_1+1-x,N-x,\bm{\lambda}_2)\n
  &\qquad\quad-\beta_1\beta_3\!\!\sum_{z_2=0}^{\min(x-1,y-2)}\!\!\!\!\!\!
  \pi(z_2,x-1,\bm{\lambda}'_1)\!\!\!\!\sum_{z_1=\max(x-2,y-2)}^{N-2}
  \!\!\!\!\!\!\!\!\!\!\pi(y-2-z_2,z_1-z_2,\bm{\lambda}'_3)\n
  &\hspace{79mm}\times
  \pi(z_1+2-x,N-x,\bm{\lambda}_2)\Bigr),
  \label{Keta4}
\end{align}
(\romannumeral5) $\eta(x)=x$ :
\begin{align}
  &\quad K(y,x)\eta(y)\n
  &=\beta_2\beta_3\eta(N)\!\sum_{z_2=0}^{\min(x,y-1)}\!\!\!\!
  \pi(z_2,x,\bm{\lambda}_1)\!\sum_{z_1=x+y-z_2}^N\!\!\!\!
  \pi(y-z_2-1,z_1-x-1,\bm{\lambda}'_3)\n
  &\hspace{65mm}\times
  \pi(z_1-x\!-1,N-x-\!1,\bm{\lambda}'_2)\n
  &\quad+\eta(x)\Bigl(
  -\beta_2\beta_3\!\sum_{z_2=0}^{\min(x,y-1)}\!\!\!\!
  \pi(z_2,x,\bm{\lambda}_1)\!\sum_{z_1=x+y-z_2}^N\!\!\!\!
  \pi(y-z_2-1,z_1-x-1,\bm{\lambda}'_3)\n
  &\hspace{72mm}\times
  \pi(z_1-x\!-1,N-x-\!1,\bm{\lambda}'_2)
  \label{Keta5}\\
  &\qquad\quad+\beta_1\sum_{z_2=1}^{\min(x,y)}\!\!
  \pi(z_2-1,x-1,\bm{\lambda}'_1)\!\sum_{z_1=x+y-z_2}^N\!\!\!\!
  \pi(y-z_2,z_1-x,\bm{\lambda}_3)\pi(z_1-x,N-x,\bm{\lambda}_2)\Bigr).
  \nonumber
\end{align}

Since each summation in \eqref{Ketas}--\eqref{Keta5} has the same structure
as the original $K(y,x)$ with shifted arguments and parameters,
the next step ($K(y,x)\eta(y)^2$) and further steps ($K(y,x)\eta(y)^3,\ldots$)
go almost parallel with the help of formulas, like
\begin{align*}
  \eta(x)=x&:\ \ \eta(y)=\eta(y-i)+\eta(i),\ \ \eta(x-j)=\eta(x)+\eta(-j),\\
  \eta(x)=q^{-x}-1&:\ \ \eta(y)=q^{-i}\eta(y-i)+\eta(i),
  \ \ \eta(x-j)=q^j\eta(x)+\eta(-j),
\end{align*}
and summation over $y$ gives \eqref{triang1}.

For type (\romannumeral2) and (\romannumeral5) with $\eta(x)=q^{-x}-1$,
we obtain $\sum_{y\in\cX}K(y,x)\eta(y)=\alpha_1+\alpha_2\eta(x)+\alpha_3
\eta(N-x)$ ($\alpha_i$ : constant), and the triangularity \eqref{triang1}
does not hold for these cases.

For the semi-infinite Markov chains given in \S\,\ref{sec:exa} and the examples
with $\eta(x)=1-q^x$ in \S\,\ref{sec:lqJ}, similar proof of triangularity holds.
Obtaining the explicit form of the coefficients $\{a_{n\,m}\}$ in
\eqref{triang1} is not necessary.
One only has to convince oneself that the triangularity holds.
The eigenvalues are easily obtained by the formula \eqref{eigform} in
{\bf Theorem \ref{theo2}}.
Since triangularity is the consequence of \eqref{pieta1}--\eqref{pieta3},
it is quite natural to expect that it also holds for convolutions other than
type (\romannumeral1)--(\romannumeral5).



\begin{thebibliography}{99}

\bibitem{nikiforov}
A.F.\,Nikiforov, S.K.\,Suslov and V.B.\,Uvarov,
{\it Classical Orthogonal Polynomials of a Discrete Variable\/},
Springer-Verlag, Berlin, (1991).

\bibitem{askey}
G.\,E.\,Andrews, R.\,Askey and R.\,Roy,
{\it Special Functions},
Encyclopedia of mathematics and its applications,
Cambridge Univ. Press, Cambridge, (1999).

\bibitem{ismail}
M.\,E.\,H.\,Ismail,
{\it Classical and Quantum Orthogonal Polynomials in One Variable\/},
Encyclopedia of mathematics and its applications,
Cambridge Univ. Press, Cambridge, (2005).

\bibitem{kls}
R.\,Koekoek, P.\,A.\,Lesky and R.\,F.\,Swarttouw,
{\it Hypergeometric orthogonal polynomials and their $q$-analogues,\/}
Springer Monographs in Mathematics,
Springer-Verlag Berlin-Heidelberg, (2010).

\bibitem{gasper}
G.\,Gasper and M.\,Rahman,
{\it Basic hypergeometric series\/}, 2nd ed.,
Encyclopedia of mathematics and its applications,
Cambridge Univ. Press, Cambridge, (2004).

\bibitem{bdsol}
R.\,Sasaki,
``Exactly Solvable Birth and Death Processes,''
J. Math. Phys. {\bf 50} (2009) 103509 (18 pp),
{\tt arXiv:0903.3097[math-ph]}.

\bibitem{os34}
S.\,Odake and R.\,Sasaki,
``Orthogonal Polynomials from Hermitian Matrices \II,''
J. Math. Phys. {\bf 59} (2018) 013504 (42pp),
{\tt arXiv:1604.00714[math.CA]}.

\bibitem{dtbd}
R.\,Sasaki,
``Exactly solvable discrete time Birth and Death processes,''
J. Math. Phys. {\bf 63} (2022) in press,
{\tt arXiv:2106.\hspace{0pt}03284[math.PR]}.

\bibitem{os12}
S.\,Odake and R.\,Sasaki,
``Orthogonal Polynomials from Hermitian Matrices,''
J. Math. Phys. {\bf 49} (2008) 053503 (43 pp),
{\tt arXiv:0712.4106[math.CA]}.
(For the dual polynomials in \S\,5.2.4 and \S\,5.2.8, see \cite{os34}.)

\bibitem{atakishi1}
M.\,N.\,Atakishiyev, N.\,M.\,Atakishiyev and A.\,U.\,Klimyk,
``Big $q$-Laguerre and $q$-Meixner polynomials and representation of the
quantum algebra $U_q(su_{1,1}$),''
J. Phys. {\bf A36} (2003) 10335-10347,
{\tt arXiv:math/0306201[math.QA]}.

\bibitem{hoa-rah83}
M.\,R.\,Hoare and M.\,Rahman,
``Cumulative Bernoulli trials and Krawtchouk processes,''
Stochastic Processes and their applications {\bf 16} (1984) 113-139.

\bibitem{coo-hoa-rah77}
R.\,D.\,Cooper, M.\,R.\,Hoare and M.\,Rahman,
``Stochastic Processes and Special Functions:
On the Probabilistic Origin of Some Positive Kernels Associated
with Classical Orthogonal Polynomials,''
J. Math. Anal. Appl. {\bf 61} (1977) 262-291.

\bibitem{albert}
F.\,A.\,Gr\"unbaum and M.\,Rahman,
``A System of Multivariable Krawtchouk Polynomials and a Probabilistic
Application,''
SIGMA {\bf 7} (2011) 119,
{\tt arXiv:1106.1835[math.PR]}.

\bibitem{diaconis20}
P.\,Diaconis and C.\,Zhong,
``Hahn polynomials and Burnside process,''
{\tt arXiv:2012.\hspace{0pt}13829[math.PR]}.

\bibitem{feller}
W.\,Feller,
``The birth and death processes as diffusion processes,''
J. Math. Pures Appl. (9) {\bf 38} (1959) 301--345.

\end{thebibliography}
\end{document}